\documentclass[sn-aps]{sn-jnl}

\usepackage{graphicx}%
\usepackage{multirow}%
\usepackage{amsmath,amssymb,amsfonts}%
\usepackage{amsthm}%
\usepackage{mathrsfs}%
\usepackage[title]{appendix}%
\usepackage{xcolor}%
\usepackage{textcomp}%
\usepackage{manyfoot}%
\usepackage{booktabs}%
\usepackage[final]{changes}

\theoremstyle{thmstyleone}%
\newtheorem{theorem}{Theorem}
\newtheorem{proposition}[theorem]{Proposition}%

\theoremstyle{thmstyletwo}%
\newtheorem{remark}{Remark}%

\theoremstyle{thmstylethree}%

\newcommand{\Cc}{\mathbb{C}}
\newcommand{\dd}{\mathrm{d}}
\newcommand{\epst}{\tilde{\epsilon}}
\newcommand{\geqs}{\geqslant}
\newcommand{\leqs}{\leqslant}
\newcommand{\Rr}{\mathbb{R}}
\newcommand{\yt}{\tilde{y}}

\newcommand{\Asymm}{A_{\mathrm{symm}}}
\newcommand{\Askew}{A_{\mathrm{skew}}}
\newcommand{\rlin}{r^{\mathsf{lin}}}

\begin{document}

\title[Convergence of waveform relaxation for nonlinear ODE systems]{On convergence of waveform relaxation for nonlinear systems of 
ordinary differential equations}

\subtitle{{\normalsize To the memory of Willem Hundsdorfer}}

\author*[1]{\fnm{M. A.} \sur{Botchev}}\email{botchev@kiam.ru}

\affil*[1]{%
\orgname{Keldysh Institute of Applied Mathematics of Russian Academy of Sciences},
\orgaddress{\street{Miusskaya Sq., 4}, 
\city{Moscow}, \postcode{125047}, \country{Russia}}}

\date{April 4, 2024}
	
\abstract{%
To integrate large systems of nonlinear differential equations in time, 
we consider a variant of nonlinear waveform relaxation 
(also known as dynamic iteration or Picard--Lindel\"of iteration), 
where at each iteration a linear inhomogeneous system of differential equations 
has to be solved. 
This is done by the exponential block Krylov 
subspace (EBK) method.
Thus, we have an inner-outer iterative method, where iterative approximations 
are determined over a certain time interval, with no time stepping involved.
This approach has recently been shown to be efficient as a time-parallel integrator
within the PARAEXP framework.
In this paper, convergence behavior of this method is assessed 
theoretically and practically.
We examine efficiency of the method by testing it on
nonlinear Burgers, 
three-dimensional Liouville--Bratu--Gelfand,
and three-dimensional nonlinear heat conduction equations 
and comparing its performance 
with that of conventional time-stepping integrators.}


\maketitle
	
\section{Introduction} 
Large systems of time-dependent ordinary differential equations arise
in various applications and in many cases have to be integrated in time
by implicit methods, see e.g.~\cite{HundsdorferVerwer:book,ElmanSilvesterWathen:book}.
Last decennia the niche of implicit methods has been gradually taken up by
exponential time integrators~\cite{HochbruckOstermann2010}.
For implicit and exponential methods the key issue is how to solve
arising nonlinear systems and (or) to evaluate related matrix functions
efficiently.  To achieve efficiency, different approaches exist and
widely applied, such as inexact Newton methods combined with powerful
linear preconditioned solvers~\cite{BrownSaad,ChoquetErhel,TromeurdervoutVassilevski2006},
splitting and Rosenbrock 
methods~\cite{Yanenko,Dyakonov64,CsomosFaragoHavasi05,ros2,HundsdorferVerwer:book} and
approximate iterative implicit schemes (which can be seen as stabilized
explicit schemes)~\cite{RKC,LokLokDAN,RKC97,TalEzer89,Lebedev98,RKC2004,MRAIpap,MRAIpar,Zhukov2011}.

Another important approach to achieve efficiency in implicit and exponential
time integrators is based on waveform relaxation methods~\cite{Lelarasmee_ea1982,NewtonSangiovanni1983,Vandewalle1993corrected},
where iterative approximations are time dependent functions
rather than time step values of a numerical solution.
These methods have also been known as dynamic iteration or
Picard--Lindel\"of iteration~\cite{MiekkalaNevanlinna1996}.
They have been developed since the 
80s~\cite{NewtonSangiovanni1983,White_ea1985,MiekkalaNevanlinna1987,LubO87,JanssenVandewalle1996}
and have received attention
recently in connection to the time-parallel exponential method PARAEXP~\cite{PARAEXP} 
and its extension to nonlinear problems~\cite{Kooij_ea2017,Kooij_ea2018}.

Typically matrix function evaluations within exponential integrators are
carried out by special linear algebra iterative procedures (often these are
Krylov subspace or Chebyshev polynomial 
methods)~\cite{ParkLight86,Henk:f(A),DruskinKnizh89,TalEzer89,GallSaad92,DruskinKnizh95,HochLub97}.
The key attractive feature of the waveform relaxation methods is that
they employ this linear algebra machinery across a certain time interval
rather than within a time step, so that computational costs are distributed
across time.  This leads to a higher computational efficiency
as well as to a parallelism across time~\cite{Vandewalle1993corrected}.

One promising time-integration exponential method of this class is
the nonlinear EBK (exponential block Krylov) method~\cite{Kooij_ea2017,Kooij_ea2018}.
It is based on block Krylov subspace projections~\cite{Botchev2013} combined with
waveform relaxation iteration employed to handle nonlinearities.
The nonlinear EBK method possesses an across-time parallelism,
and, in its parallel implementation, can be seen as an efficient version
of the time-parallel PARAEXP method~\cite{Kooij_ea2017,Kooij_ea2018}.

Convergence of the waveform relaxation methods has been studied
in different settings, in particular, for linear initial-value problems,
for nonlinear Gauss--Seidel and Jacobi iterations (the working horses of
classical waveform relaxation) and
for time-discretized settings, see~\cite{MiekkalaNevanlinna1996,Vandewalle1993corrected} for a survey.
Convergence results for waveform relaxation in general nonlinear settings
are scarce, and it is often assumed that
\emph{``studying the linear case carefully \dots{}
might be what users really need''}~\cite{MiekkalaNevanlinna1996}.
Except book~\cite{Vandewalle1993corrected}, containing some
general introductory convergence results, 
one of the papers where nonlinear waveform relaxation is
studied is~\cite{NevanlinnaOdeh1987}.  Below, in Remark~\ref{NevOdeh}, 
we comment more specifically on the results presented there.
The aim of this paper is to extend convergence results for waveform relaxation 
to a nonlinear setting employed in the EBK method.
This also yields an insight on its convergence behavior in practice.

The paper is organized as follows.
In Section~\ref{s:setting} a problem setting and the assumptions made are
formulated.
Main convergence results are presented in Section~\ref{s:conv}.
Since the results are formulated for the error function
which is in general unknown, in Section~\ref{s:nonl_r} an error estimate
in terms of a computable nonlinear residual is formulated.
At each waveform relaxation iteration a linear initial-value problem (IVP)
has to be solved which, in current setting, is done by an iterative block
Krylov subspace method (the linear EBK method).
Therefore, in Section~\ref{s:lin_r} we show how the nonlinear iteration
convergence is affected if the linear IVP is solved inexactly.
In Section~\ref{s:implem} implementation issues are briefly discussed.
Numerical experiments are described in Section~\ref{s:num_exp1D}
(for a 1D Burgers equation), \deleted{and}
Section~\ref{s:num_exp3D} (for a 3D Liouville–Bratu–Gelfand equation),
\added{and Section~\ref{s:num_heat} 
(for a 3D nonlinear heat conduction problem)}.
We finish the paper with conclusions drawn in Section~\ref{s:concl}.

\section{Nonlinear waveform relaxation and its convergence}
\subsection{Problem setting, assumptions and notation}
\label{s:setting}
We consider an IVP for nonlinear ODE system
\begin{equation}
\label{IVP0}
y'(t) = \Phi(t,y(t)),  \quad y(0)=v,\quad t\in[0,T]
\end{equation}
where $\Phi:\Rr\times\Rr^N\rightarrow\Rr^N$, $v\in\Rr^N$ and $T>0$ are given.
Without loss of generality, we assume that
the time dependence in the right hand side function $\Phi$ is of the form
\begin{equation}
\label{autonom}
\Phi(t,y) = \bar{\Phi}(y) + g(t),
\quad \bar{\Phi}:\Rr^N\rightarrow\Rr^N,
\quad g : \Rr\rightarrow\Rr^N.
\end{equation}
Otherwise, we can transform~\eqref{IVP0} to an equivalent autonomous form
(usually this is done by extending the unknown vector function $y(t)$ with an additional
$(N+1)$-th coordinate $y_{N+1}(t)=t$ and adding an ODE $y_{N+1}'(t)=1$ to the
system).

Iterative methods considered in this note are based on a family
of splittings
\begin{equation}
\label{split_f}
\Phi(t,y) = -A_k y + f_k(y) + g(t), \qquad \forall\; y\in\Rr^N, \; t\geqs 0,
\end{equation}
with $k$ being the iteration number, so that matrices $A_k\in\Rr^{N\times N}$
and mappings $f_k:\Rr^N\rightarrow\Rr^N$
may vary from iteration to iteration.
The mappings $f_k$ in~\eqref{split_f} are all assumed to be Lipschitz continuous with 
a common Lipschitz constant $L$, i.e.,
\begin{equation}
\label{L}
\forall k:\quad
\|f_k(u) - f_k(v)\| \leqs L \|u-v\|, \qquad \forall\; u,v\in\Rr^N.
\end{equation}
Here and throughout the paper $\|\cdot\|$ denotes a vector norm in $\Cc^N$
or a corresponding induced matrix norm.
Furthermore, we assume that for all the matrices $A_k$ in~\eqref{split_f}
there exist constants $C>0$ and $\omega\geqs 0$ such that 
\begin{equation}
\label{A}
\|\exp(-tA_k)\|\leqs C e^{-\omega t}, \qquad t\geqs 0.
\end{equation}
For a particular matrix $A=A_k$, condition~\eqref{A} holds, for instance, in the 2-norm 
if the numerical range of $A$ lies in the complex halfplane 
$\{z=x+iy \;|\; x\geqs 0, y\in\Rr, i^2=-1\}$.
In this case $C=1$ and $\omega$ is the smallest eigenvalue
of the symmetric part $\frac12(A+A^T)$ of $A$.
Condition~\eqref{A} is closely related to the concept of the matrix
logarithmic norm, see, e.g.,~\cite{Dekker-Verwer:1984,HundsdorferVerwer:book}.
In particular, if $\mu(-A)$ is the logarithmic norm of $-A$ 
then
condition~\eqref{A} holds, with $C=1$ and for all $t\geqs 0$, 
if and only if \cite[Theorem~I.2.4]{HundsdorferVerwer:book} 
$$
\mu(-A)\leqs -\omega .
$$
In the analysis below, we also use functions 
$\varphi_j$~\cite[relation~(2.10)]{HochbruckOstermann2010}, 
defined for $j=0,1,2,\dots$ as
\begin{equation}
\label{phi}
\varphi_{j+1}(z)=\frac{\varphi_j(z)-\varphi_j(0)}{z},
\qquad \varphi_0(z)=e^z,  
\end{equation}
where we set $\varphi_j(0)=1/j!$, which makes $\varphi_j$ entire
functions.
In this paper, the following variation-of-constants formula (see,
e.g.,~\cite[Section~I.2.3]{HundsdorferVerwer:book},
\cite[relation~(1.5)]{HochbruckOstermann2010}) is instrumental:
\begin{equation}
\label{VOC}
y(t) = \exp(-tA_k)v + \int_0^t \exp(-(t-s)A_k)
\left[f_k(y(s)) + g(s) \right] \dd s,
\quad t\in[0,T],
\end{equation}
where $y(t)$ is solution of IVP~\eqref{IVP0},\eqref{split_f}.

\subsection{Nonlinear waveform relaxation}
\label{s:conv}
Let $y_0(t)$ be an approximation to the unknown function $y(t)$
for $t\in[0,T]$, with $y_0(0)=v$.
Usually $y_0(t)$ is taken to be $y_0(t)\equiv v$ for all $t\in[0,T]$.
To solve~\eqref{IVP0}, in this paper we consider nonlinear waveform relaxation iteration
where we solve successfully, 
for $k=0,1,2,\dots$, a linear inhomogeneous IVP
\begin{equation}
\label{nonlR}
y_{k+1}'(t) = -A_k y_{k+1}(t) + f_k(y_k(t)) + g(t),  
\quad y_{k+1}(0)=v,\quad t\in[0,T].
\end{equation}
Here the matrices $A_k\in\Rr^{N\times N}$ and the mappings $f_k:\Rr^N\rightarrow\Rr^N$ 
form a splitting of $\Phi$, see~\eqref{split_f}, and satisfy, for all $k$,
the assumptions~\eqref{L},\eqref{A}.
We emphasize that at each iteration an approximation $y_{k+1}(t)$ is computed
and stored for the whole time range $t\in[0,T]$.  In Section~\ref{s:implem}
below we briefly discuss how to do this
efficiently. 

The following proposition provides a sufficient condition for 
iteration~\eqref{nonlR} to converge.  

\begin{proposition}
\label{prop1}
Let IVP~\eqref{IVP0} be solved iteratively by~\eqref{nonlR} and let 
assumptions~\eqref{L},\eqref{A} hold for $A_k$ and $f_k$ in~\eqref{nonlR}, $k=0,1,\dots$.
Then for the error $\epsilon_{k+1}(t)\equiv y(t) - y_{k+1}(t)$
of the iterative approximation $y_{k+1}(t)$ holds, for $k=0,1,\dots$,
\begin{equation}
\label{est1}
\|\epsilon_{k+1}(t)\|\leqs C L t\varphi_1(-\omega t)
\max_{s\in[0,t]} \|\epsilon_k(s)\|, \quad \forall t\in[0,T],
\end{equation}
and nonlinear waveform relaxation~\eqref{nonlR} converges for $t\in[0,T]$
to solution $y(t)$ of IVP~\eqref{IVP0},\eqref{split_f} provided that
\begin{equation}
\label{conv}
C L T\varphi_1(-\omega T) <1.  
\end{equation}
\end{proposition}

\begin{proof}
Subtracting the iteration formula~\eqref{nonlR}
from the ODE $y'(t)=-A_k y(t)  + f_k(y(t)) + g(t)$ and taking into
account that $y(0)=y_{k+1}(0)=v$, we see that the error 
$\epsilon_{k+1}(t)$ satisfies IVP
\begin{equation}
\label{ivp_err}
\epsilon_{k+1}'(t) = -A_k \epsilon_{k+1}(t) + f_k(y(t)) - f_k(y_k(t)),
\quad
\epsilon_{k+1}(0) = 0.
\end{equation}
Then, applying the variation-of-constants formula~\eqref{VOC}
to IVP~\eqref{ivp_err}, we obtain
$$
\epsilon_{k+1}(t) = \int_0^t \exp(-(t-s)A_k)\left[f_k(y(s))-f_k(y_k(s))\right]\dd s.
$$
Therefore, using~\eqref{L} and~\eqref{A}, we can bound
$$
\begin{aligned}
\|\epsilon_{k+1}(t)\| 
&\leqs 
\int_0^t \| \exp(-(t-s)A_k)\|\,\|\left[f_k(y(s))-f_k(y_k(s))\right]\|\dd s
\\
&\leqs 
\int_0^t \| \exp(-(t-s)A_k)\|\,L \|y(s))-y_k(s)\|\dd s
\\
&\leqs 
L\max_{s\in[0,t]}\|\epsilon_k(s)\| \int_0^t \| \exp(-(t-s)A_k)\|\,\dd s
\\
&\leqs C L\max_{s\in[0,t]}\|\epsilon_k(s)\| \int_0^t e^{-\omega(t-s)}\,\dd s
= C L t \varphi_1(-\omega t) \max_{s\in[0,t]}\|\epsilon_k(s)\|.
\end{aligned}
$$
Thus, \eqref{est1} is proved.  Taking into account that $t \varphi_1(-\omega t)$
is a monotonically increasing function of $t$, we obtain~\eqref{conv}.
\end{proof}

\begin{remark}
Proposition~\ref{prop1} shows that we can choose a length~$T$ of the time interval $[0,T]$
such that the nonlinear waveform relaxation converges.
The larger the Lipschitz constant $L$, the smaller $T$ should be taken.
\end{remark}

\begin{remark}
\label{NevOdeh}  
To solve an initial-value problem for an autonomous ODE system
$$
C(y)y'(t) -\bar{F}(y(t)) = 0,
$$
in~\cite{NevanlinnaOdeh1987} a nonlinear waveform iteration 
\begin{equation}
\label{NevOdeh_it}
y'_{k+1}(t) - \bar{F}(y_{k+1}(t)) = \widehat{g}(y'_k(t))
\end{equation}
is considered, where $C(y)$ is a matrix and $C(y)y' = y' - \widehat{g}(y')$.
A convergence analysis of~\eqref{NevOdeh_it} is given in assumption that 
$\widehat{g}$ is Lipschitz-continuous with a Lipschitz constant less than one.
In~\cite{NevanlinnaOdeh1987}, a particular case of~\eqref{NevOdeh_it}
is also considered with $\bar{F}(y) = -Ay + f(y)$, i.e.,
$$
y'_{k+1}(t) + Ay_{k+1}(t) - f(y_{k+1}(t)) = \widehat{g}(y'_k(t)).
$$
Hence, our results here do not overlap with the results in~\cite{NevanlinnaOdeh1987}. 
\end{remark}

\medskip

In case $f(y)$ is linear, waveform relaxation~\eqref{nonlR} is known
to have an attractive property to converge superlinearly for 
final~$T$, see, e.g.,~\cite[Theorem~5.1]{BGH13} 
and, for inexact waveform relaxation, \cite[Theorem~2]{BotchevOseledetsTyrtyshnikov2014}.
The following result shows that the superlinear convergence property 
is shared by nonlinear waveform relaxation~\eqref{nonlR}.

\begin{proposition}
\label{prop2}
Let IVP~\eqref{IVP0} be solved iteratively by~\eqref{nonlR} and let 
assumptions~\eqref{L},\eqref{A} hold for $A_k$ and $f_k$ in~\eqref{nonlR}, $k=0,1,\dots$.
Then for the error $\epsilon_k(t)$ 
of the iterative approximation $y_k(t)$ holds, for $k=1,2,\dots$,
\begin{equation}
\label{est2}
\begin{aligned}
\|\epsilon_k(t)\| &\leqs (C L)^k t^k e^{-\omega t} \varphi_k(\omega t)
\max_{s\in[0,t]} \|\epsilon_0(s)\|
\\
&\leqs (C L)^k \frac{t^k}{k !} 
\max_{s\in[0,t]} \|\epsilon_0(s)\|, \quad \forall t\in[0,T].
\end{aligned}
\end{equation}
\end{proposition}
 
\begin{proof}
The proof is very similar to the proof of~\cite[Theorem~5.1]{BGH13}, where 
a linear equivalent of~\eqref{IVP0},\eqref{nonlR} is considered.
The estimate~\eqref{est2} will be proven by induction.
Note that $t\varphi_1(-\omega t)= te^{-\omega t}\varphi_1(\omega t)$.
Therefore, the estimate~\eqref{est2} for $k=1$ holds, as it coincides 
with~\eqref{est1} for $k=0$.
Assume that~\eqref{est2} holds for a certain~$k$.
From the proof of Proposition~\ref{prop1}, we see that
\begin{align*}
\|\epsilon_{k+1}(t)\| &\leqs C L \int_0^t e^{-(t-s)\omega}\|\epsilon_{k}(s)\|\dd s
\\
\intertext{and, by the induction assumption,}
&\leqs
C L \int_0^t e^{-(t-s)\omega}(C L)^k s^k e^{-\omega s} \varphi_k(\omega s)
\max_{\tilde{s}\in[0,s]} \|\epsilon_0(\tilde{s})\|\dd s
\\
&=
(C L)^{k+1} e^{-\omega t} \int_0^t s^k \varphi_k(\omega s)
\max_{\tilde{s}\in[0,s]} \|\epsilon_0(\tilde{s})\|\dd s
\\
&\leqs
(C L)^{k+1} e^{-\omega t} \max_{s\in[0,t]} \|\epsilon_0(s)\|
\int_0^t s^k \varphi_k(\omega s)\dd s 
\\
&=
(C L)^{k+1} t^{k+1} e^{-\omega t} \varphi_{k+1}(\omega t) 
\max_{s\in[0,t]} \|\epsilon_0(s)\|,
\end{align*}
where the relation 
$\int_0^t s^k \varphi_k(\omega s)\dd s = t^{k+1} \varphi_{k+1}(\omega t)$
is employed.
Thus, the induction step is done.
Finally, using the definition of $\varphi_k$ and the fact
that $\omega\geqs 0$, it is not difficult to check that 
$$
t^k e^{-\omega t}\varphi_k(\omega t) \leqs \frac{t^k}{k!},
\quad t\geqs 0,
$$
which proves the second inequality in~\eqref{est2}.
\end{proof}

\subsection{Residual of the nonlinear problem}
\label{s:nonl_r}
To control convergence of waveform relaxation~\eqref{nonlR}
in practice, it is natural to consider a residual $r_k(t)$ of an iterative
approximation $y_k(t)$ in~\eqref{nonlR}.
Since $y_k(t)$ is defined by~\eqref{nonlR},\eqref{split_f} written for $k-1$,
it is natural to define the residual with respect to IVP~\eqref{IVP0}
combined with splitting~\eqref{split_f} for $k-1$, i.e., for IVP
\begin{equation}
\label{IVPk}  
y'(t) = -A_{k-1}y(t) + f_{k-1}(y(t)) + g(t),\quad
y(0) = v, \quad t\in[0,T].
\end{equation}
Hence, we define
\begin{equation}
\label{rk_def}
r_k(t) \equiv -A_{k-1}y_k(t) + f_{k-1}(y_k(t)) + g(t) - y_k'(t), \quad t\in [0,T],
\end{equation}
and write
\begin{equation}
\label{rk}
\begin{aligned}
r_k(t) &= -A_{k-1}y_k(t) + f_{k-1}(y_k(t)) + g(t) - y_k'(t) \pm f_{k-1}(y_{k-1}(t)) \\
       &= \underbrace{-A_{k-1}y_k(t) + f_{k-1}(y_{k-1}(t)) + g(t) - y_k'(t)}_{= 0} + 
          f_{k-1}(y_k(t)) - f_{k-1}(y_{k-1}(t)) \\
       &= f_{k-1}(y_k(t)) - f_{k-1}(y_{k-1}(t)), \quad t\in[0,T], \quad k=1,2,\dots .
\end{aligned}
\end{equation}
Note that the residual $r_k(t)$ possesses the backward error property,
i.e., iterative approximation $y_k(t)$ can be viewed as the exact solution of
a perturbed IVP~\eqref{IVPk}
\begin{equation}
\label{IVP_rk}
y_k'(t) = -A_{k-1}y_k(t) + f_{k-1}(y_k(t)) + g(t) - r_k(t),
\quad y_k(0) = v \quad t\in [0,T].  
\end{equation}
Subtracting the ODE of this problem from the ODE in~\eqref{IVPk}, we obtain
an IVP for the error $\epsilon_k(t)$
\begin{equation}
\label{ivp_err1}
\epsilon_k'(t) = -A_{k-1}\epsilon_k(t) + f_{k-1}(y(t))-f_{k-1}(y_k(t)) + r_k(t),
\quad
\epsilon_k(0) = 0,
\end{equation}
which is equivalent to the IVP~\eqref{ivp_err}.
The following proposition shows that the residual can be seen 
as an upper bound bound for the error.

\begin{proposition}
\label{r2e}
Let IVP~\eqref{IVP0} be solved iteratively by~\eqref{nonlR} and let 
assumptions~\eqref{L},\eqref{A} hold for $A_k$ and $f_k$ in~\eqref{nonlR}, $k=0,1,\dots$.
Let $T>0$ be chosen such that the sufficient condition~\eqref{conv} for 
convergence of~\eqref{nonlR} holds
\begin{equation}
\label{conv1}
C L T\varphi_1(-\omega T) \leqs \delta <1,
\end{equation}
for a certain constant $\delta\in(0,1)$.
Then, for $\forall t\in[0,T]$,
\begin{equation}
\label{r2e:eq}
\max_{s\in[0,t]} \|\epsilon_k(s)\| \leqs 
\frac{C t\varphi_1(-\omega t)}{1-C L t\varphi_1(-\omega t)} \max_{s\in[0,t]}\|r_k(s)\|
\leqs 
\frac{\delta}{(1- \delta)L} \max_{s\in[0,t]}\|r_k(s)\|.
\end{equation}
Note that $T\varphi_1(-\omega T)$ increases with $T$ monotonically.
\end{proposition}

\begin{proof}
Employing the variation-of-constants formula for IVP~\eqref{ivp_err1},
we can bound
\begin{align*}
\|\epsilon_k(t)\| &\leqs 
\int_0^t \left\| \exp(-(t-s)A_{k-1})\left[f_{k-1}(y(s))-f_{k-1}(y_k(s)) + r_k(s)\right] \right\|\dd s
\\&\leqs 
\max_{s\in[0,t]}\left\|f_{k-1}(y(s))-f_{k-1}(y_k(s)) + r_k(s)\right\| \int_0^t \| \exp(-(t-s)A_{k-1})\|\,\dd s  
\\&\leqs 
\max_{s\in[0,t]}\left\|f_{k-1}(y(s))-f_{k-1}(y_k(s)) + r_k(s)\right\| \, C t\varphi_1(-\omega t)
\\&\leqs 
\left( L\max_{s\in[0,t]}\|\epsilon_k(s)\| + \max_{s\in[0,t]}\|r_k(s)\|\right) 
\, C t\varphi_1(-\omega t)
\\&\leqs 
C L t\varphi_1(-\omega t) \max_{s\in[0,t]}\|\epsilon_k(s)\| + 
C t\varphi_1(-\omega t) \max_{s\in[0,t]}\|r_k(s)\|. 
\end{align*}
Taking into account that 
$C L t\varphi_1(-\omega t)\leqs C L T\varphi_1(-\omega T) \leqs \delta <1$
for $t\in[0,T]$, 
we obtain
\begin{align*}
(1-C L t\varphi_1(-\omega t))\max_{s\in[0,t]} \|\epsilon_k(s)\| &\leqs 
C t\varphi_1(-\omega t) \max_{s\in[0,t]}\|r_k(s)\|
\\
\max_{s\in[0,t]} \|\epsilon_k(s)\| &\leqs 
\frac{C t\varphi_1(-\omega t)}{1-C L t\varphi_1(-\omega t)} \max_{s\in[0,t]}\|r_k(s)\|,
\end{align*}
which yields~\eqref{r2e:eq}.
\end{proof}

\added{As stated by proposition below, an upper bound for the nonlinear residual can 
be obtained in terms of the error at the previous iteration.}

\begin{proposition}
\label{e2r}
\added{Let IVP~\eqref{IVP0} be solved iteratively by~\eqref{nonlR} and let 
assumptions~\eqref{L},\eqref{A} hold for $A_k$ and $f_k$ in~\eqref{nonlR}, $k=0,1,\dots$.
Then, for $k=1,2,\dots$,
\begin{align}
\label{e2r:eq1}
\|r_k(t)\| &\leqs (1 + C L t\varphi_1(-\omega t) )
L \max_{s\in[0,t]} \|\epsilon_{k-1}(s)\|,
\quad \forall t\in[0,T],
\\
\intertext{If, in addition, convergence condition~\eqref{conv1} is satisfied,
we have}
\label{e2r:eq2}
\|r_k(t)\| &\leqs (1 + \delta ) L \max_{s\in[0,t]} \|\epsilon_{k-1}(s)\|
< 2 L \max_{s\in[0,t]} \|\epsilon_{k-1}(s)\|, \quad \forall t\in[0,T].
\end{align}}
\end{proposition}

\begin{proof}
\added{Using relation~\eqref{rk}, we can write, for $t\in[0,T]$,
\begin{align*}
r_k(t) &= f_{k-1}(y_k(t)) - f_{k-1}(y(t)) + f_{k-1}(y(t)) - f_{k-1}(y_{k-1}(t)),
\\
r_k(t) &\leqs \|f_{k-1}(y_k(t)) - f_{k-1}(y(t))\| + \|f_{k-1}(y(t)) - f_{k-1}(y_{k-1}(t))\|
\\
&\leqs L\|y_k(t)) - y(t)\| + L\|y(t) - y_{k-1}(t)\|
= L(\|e_k(t)\| + \|e_{k-1}(t)\|)
\\
&\leqs L\left(
CL t\varphi_1(-\omega t)\max_{s\in[0,t]}\|e_{k-1}(s)\| + \|e_{k-1}(t)\|
\right),
\end{align*}
which leads to~\eqref{e2r:eq1},\eqref{e2r:eq2}.}
\end{proof}

\subsection{Linear inner iteration and its residual}
\label{s:lin_r}
In practice, the linear ODE system~\eqref{nonlR} at each waveform relaxation 
iteration can be solved by any suitable integrator 
\added{inexactly,} with a certain accuracy tolerance.  
In the context of the time-parallel PARAEXP scheme, 
the nonlinear waveform relaxation
appeared to be efficient for fluid dynamics problems~\cite{Kooij_ea2017,Kooij_ea2018},
with linear IVP~\eqref{nonlR} solved 
by the exponential block Krylov subspace (EBK) method~\cite{Botchev2013}. 
Thus, in this setting we have an inner-outer iterative process where 
each nonlinear waveform relaxation iteration~\eqref{nonlR} involves an inner Krylov
subspace iterative process $\tilde{y}_{k+1,\ell}(t)\rightarrow y_{k+1}(t)$,
with $\ell$ being the inner iteration number.
For notation simplicity we omit the dependence on the inner iteration number $\ell$
in $\tilde{y}_{k+1,\ell}(t)$ and write $\tilde{y}_{k+1}(t)$.
These inner EBK iterations are controled by checking the norm of a residual
defined, for $\tilde{y}_{k+1}(t)\approx y_{k+1}(t)$, as
\begin{equation}
\label{rm}
\rlin_{k+1}(t)\equiv -A_k\tilde{y}_{k+1}(t)+ f_k(\tilde{y}_k(t)) + g(t) 
-\tilde{y}'_{k+1}(t), \quad t\in[0,T].
\end{equation}
Here we assume that 
the previous iterative approximation \deleted{$y_k(t)$} is also computed inexactly and
therefore we replace $y_k(t)$ by $\tilde{y}_k(t)$.
The inner iterations stop after a certain number $\ell$ of inner 
iterations as soon as the linear residual $\rlin_{k+1}(t)$,  
defined by~\eqref{rm}, is small enough in norm.
The residual  $\tilde{r}_{k+1}(t)$ turns out to be easily computable
as a by-product of the block Krylov subspace method~\cite{Botchev2013}.
Moreover, the residual in~\eqref{rm} can be used as a controller 
for the error in the linear IVP~\eqref{nonlR}, 
see~\cite[relation~(16)]{Botchev2013}.  


\replaced{As Proposition~\ref{P:inex} below}{The following result} shows, \deleted{that} 
convergence in the nonlinear waveform relaxation
is retained if the IVP~\eqref{nonlR} at each nonlinear iteration is solved 
approximately, such that residual~\eqref{rm} is small enough in norm.
\added{To formulate the proposition, it is convenient to assume that 
the linear residual }
\deleted{i.e., more precisely,}\added{is bounded in norm by the error achieved at
a previous nonlinear iteration, cf.~\eqref{rm_cond}.}

\begin{proposition}
\label{P:inex}
Let IVP~\eqref{IVP0} be solved iteratively by~\eqref{nonlR} and 
let assumptions~\eqref{L},\eqref{A}
and convergence condition~\eqref{conv1} hold. 
Furthermore, let the IVP~\eqref{nonlR} at each nonlinear iteration $k = 0, 1,\dots$ 
be solved approximately with $\ell$ inner iterations such that for the inner 
iteration residual~\eqref{rm} one of the following two conditions holds:
\begin{equation}
\label{rm_cond}
\begin{aligned}
\mathrm{(a)}\qquad   
& \|\rlin_{k+1}(t)\|\leqslant\eta\max_{s\in[0,t]}\|f_k(\yt_{k+1}(s)) - f_k(\yt_k(s))\| ,
\quad t\in [0,T],
\\
\mathrm{(b)}\qquad   
& \|\rlin_{k+1}(t)\|\leqslant\eta L \max_{s\in[0,t]}\|\yt_{k+1}(s) - \yt_k(s)\| ,
\quad t\in [0,T],
\end{aligned}
\end{equation}
with a constant $\eta\in(0,1)$.  
Then for the error $\epst_{k+1}(t)\equiv y(t) - \yt_{k+1}(t)$
of the iterative solution $\yt_{k+1}(t)$ in inexact inner-outer 
iteration~\eqref{nonlR},\eqref{rm} holds, for $k=0,1,\dots$,
\begin{equation}
\label{est4}
\|\epst_{k+1}(t)\| \leqs 
\frac{(1+\eta) C L t\varphi_1(-\omega t)}{1 - \eta C L t\varphi_1(-\omega t)}
\max_{s\in[0,t]} \|\epst_k(s)\|, \quad \forall t\in[0,T].
\end{equation}
Moreover, the inexact iteration~(8),(23) 
converges for $t\in[0,T]$
to solution $y(t)$ of IVP~(1),(3) 
provided that
\begin{equation}
\label{conv4}
\eta < \frac{1-\delta}{2\delta}.
\end{equation}
\end{proposition}

\begin{proof}
Note that, due to~\eqref{L},
$$
\|f_k(\yt_{k+1}(t)) - f_k(\yt_k(t))\| \leqs L \|\yt_{k+1}(t) - \yt_k(t)\|
$$  
and, hence, condition~\eqref{rm_cond},(b)
follows from condition~\eqref{rm_cond},(a)
Therefore, without loss of generality, we assume that~\eqref{rm_cond},(a)
holds.
It follows from~\eqref{rm} 
that $\tilde{y}_{k+1}(t)$ solves a perturbed ODE system
\begin{equation}
\label{pertODE}  
\tilde{y}'_{k+1}(t)= -A_k\tilde{y}_{k+1}(t)+ f_k(\tilde{y}_k(t))  + g(t) 
- \rlin_{k+1}(t). 
\end{equation}
Subtracting this equation from the original nonlinear ODE
$y'(t)= -A_k y(t)+ f_k(y(t)) + g(t)$ 
and taking into account that $\tilde{y}_{k+1}(0)=y(0)=v$,
we obtain an IVP for the error
$\epst_{k+1}(t)\equiv y(t) - \tilde{y}_{k+1}(t)$ (cf.~\eqref{ivp_err}):
$$
\epst'_{k+1}(t) = -A_k \epst_{k+1}(t) + f_k(y(t)) -
f_k(\tilde{y}_k(t)) + \rlin_{k+1}(t),
\quad \epst_{k+1}(0)=0.
$$
Application of variation-of-constants formula~\eqref{VOC}
leads to
$$
\epst_{k+1}(t) = \int_0^t \exp(-(t-s)A_k)\left[f_k(y(s))-f_k(\tilde{y}_k(s))
  + \rlin_{k+1}(s)\right]\dd s.
$$
Then, using~\eqref{L},\eqref{A},\eqref{rm_cond} 
we get, for $t\in [0,T]$,
$$
\begin{aligned}
\|\epst_{k+1}(t)\| 
&\leqs 
\int_0^t \| \exp(-(t-s)A_k)\|\, \left\| f_k(y(s))-f_k(\tilde{y}_k(s)) + \rlin_{k+1}(s)\right\|\dd s
\\
&\leqs C t \varphi_1(-\omega t) \left[ L \max_{s\in[0,t]}\|\epst_k(s)\|
+ \max_{s\in[0,t]}\|\rlin_{k+1}(s)\|\right]
\\
&\leqs C t \varphi_1(-\omega t) \left[ L \max_{s\in[0,t]}\|\epst_k(s)\|
+ \eta L \max_{s\in[0,t]}\| \yt_{k+1}(s) - \yt_k(s) \|\right].
\end{aligned}
$$
Furthermore,
$$
\max_{s\in[0,t]}\| \yt_{k+1}(s) - \yt_k(s) \pm y(s)\|
\leqs \max_{s\in[0,t]}\| \epst_{k+1}(s)\| + \max_{s\in[0,t]}\| \epst_k(s)\|,
$$
so that, for $t\in [0,T]$,
\begin{gather*}
\|\epst_{k+1}(t)\| \leqs C L t \varphi_1(-\omega t) \left[ 
(1+\eta)\max_{s\in[0,t]}\|\epst_k(s)\| + \eta \max_{s\in[0,t]}\| \epst_{k+1}(s)\|\right],
\\
(1- C L t \varphi_1(-\omega t)\eta)\max_{s\in[0,t]}\| \epst_{k+1}(s)\|
\leqs
C L t \varphi_1(-\omega t) (1+\eta)\max_{s\in[0,t]}\|\epst_k(s)\|.
\end{gather*}
Taking into account that 
$C L t \varphi_1(-\omega t)\leqs C L T \varphi_1(-\omega T)\leqs\delta<1$ and
$0<\eta<1$, we then obtain
$$
\max_{s\in[0,t]}\| \epst_{k+1}(s)\|
\leqs
\frac{C L t \varphi_1(-\omega t) (1+\eta)}{1- C L t \varphi_1(-\omega t)\eta}
\max_{s\in[0,t]}\|\epst_k(s)\|
\leqs
\frac{\delta (1+\eta)}{1- \delta\eta}
\max_{s\in[0,t]}\|\epst_k(s)\|.
$$
This yields~\eqref{est4}.  Condition~\eqref{conv4} is equivalent
to $\frac{\delta (1+\eta)}{1- \delta\eta}<1$. 
\end{proof}

\begin{remark}
Note that for $\eta\rightarrow 0$ 
the obtained convergence condition~\eqref{conv4} transforms 
into~\eqref{conv1}.
In general, condition~\eqref{rm_cond},(b) is probably easier to check
in practice than~\eqref{rm_cond},(a). 
However, IVPs where~\eqref{rm_cond},(a)
is easy to compute may occur as well.
Since $\yt_k(t)$ and $\yt_{k+1}(t)$ are solutions of IVPs with the same 
initial vector $v$, we may expect that
$\max_{s\in[0,t]}\|\yt_{k+1}(s) - \yt_k(s)\|\approx \|\yt_{k+1}(t) - \yt_k(t)\|$.

Nevertheless, employment of condition~\eqref{rm_cond} in practice
does not seem to be necessary:
as numerical tests below demonstrate, simple stopping criteria
(cf.~\eqref{lin_stop},\eqref{lin_s2})
based on smallness of $\|\rlin_{k+1}(t)\|$ work fine in practice,
leading to an efficient inner-outer iteration method.
\end{remark}


\subsection{Implementation of nonlinear waveform relaxation}
\label{s:implem}
We take initial guess function $y_0(t)$ identically (i.e., for all $t\in[0,T]$)
equal to the initial vector $y_0(t)\equiv v$.
Then, at each iteration~\eqref{nonlR} we set $\tilde{g}(t):=f_k(y_k(t))+g(t)$ 
and solve a linear IVP
\begin{equation}
\label{nonlRa}
y'_{k+1}(t) = -A_k y_{k+1}(t) + \tilde{g}(t), \quad y_{k+1}(0)=v, 
\quad t\in[0,T].  
\end{equation}
Since this IVP is solved by the EBK method~\cite{Botchev2013},
the solution $y_{k+1}(t)$ at every iteration is obtained as a linear
combination of Krylov subspace vectors, 
\added{stored columnwise in a matrix $V_{\ell\cdot m}\in\Rr^{N\times \ell m}$}, 
with $\ell$ being the number of inner iterations,
$$
\added{y_{k+1}(t)=V_{\ell\cdot m} u_{k+1}(t), \quad t\in[0,T],
\quad u_{k+1}:\Rr\rightarrow \Rr^{\ell\cdot m}},  
$$
\added{where} $u_{k+1}(t)$ of is solution of a small-sized projected 
IVP \added{(see~\cite{Botchev2013} for details)}.
This allows to store $y_{k+1}(t)$ in a compact way for all $t\in[0,T]$.
At each nonlinear iteration~$k$ the EBK method stops as soon as
the residual~\eqref{rm} is small enough in norm.
The residual norm is easily computable and is checked for several values of $t\in[0,T]$.

Starting with the second iteration $k=1$, before the IVP~\eqref{nonlRa} is solved, 
the vector function $\tilde{g}(t)\equiv f_k(y_k(t))+g(t) :\Rr\rightarrow\Rr^N$ is sampled at $n_s$ points $t_1$, \dots, $t_{n_s}$ 
covering the time interval $[0,T]$.   Usually, it is sensible to take $t_1=0$,
$t_{n_s}=T$ and the other $t_j$ to be the Chebyshev polynomial roots: 
$$
t_j = \frac{T}{2}\left(1-\cos\frac{\pi(j-3/2)}{s-2}\right),\quad
j=2,\dots,n_s-1. 
$$
The number of sample points $n_s$ should be taken such that 
$\tilde{g}(t)$ is sufficiently well approximated by its 
linear interpolation at $t_1$, \dots, $t_{n_s}$.
The computed samples $\tilde{g}(t_j)$, $j=1,\dots,{n_s}$, 
are then stored as 
the columns of a matrix $\tilde{G}\in\Rr^{N\times n_s}$ and the thin singular
value decomposition of $\tilde{G}$ is employed to obtain a low rank
representation
\begin{equation}
\label{Upt}  
\tilde{g}(t) \approx U p(t), \quad U\in\Rr^{N\times m}, \quad t\in[0,T],
\end{equation}
where typically $m<10$.
For more details of this procedure we refer to~\cite{Botchev2013}
or to~\cite{Kooij_ea2017}.
Usually the computational work needed to obtain representation
$\tilde{g}(t) \approx U p(t)$ is negligible with respect to the other work
in the EBK method.
The nonlinear iteration~\eqref{nonlR} is stopped as soon
as the nonlinear residual~\eqref{rk} is small enough in norm
(in our limited experience, in practice it suffice to check 
the residual norm at final time $t=T$ only).
\added{The resulting algorithm is outlined below.}

\clearpage

\vspace*{2ex}
\noindent
\hrule\vspace*{0.1em}\hrule\noindent
\textbf{Algorithm~1}
\hrule
\noindent
Given: IVP~\eqref{IVP0}, initial value vector $v$, tolerance
$\texttt{tol}$, block size $m$, sample number $n_s$.
\\
Result: numerical solution $y(T)$ of IVP~\eqref{IVP0}.
\\
Initialize nonlinear iteration:\\
$y_0(t):=v$ (constant in time)
\\
for $k=0,1,2,\dots$ (main iteration loop)
\\\hspace*{1em}
if ($k=0$) then
\\\hspace*{2em}
$\texttt{resnorm} := \|r_0(T)\| = \|\Phi(T,v)\|$
\\\hspace*{1em}
else
\\\hspace*{2em}
$\texttt{resnorm} := \|r_k(T)\| = \|f_{k-1}(y_k(T))-f_{k-1}(y_{k-1}(T))\|$
\\\hspace*{1em}
endif
\\\hspace*{1em}
if \texttt{resnorm} is small enough (cf.~\eqref{stop},\eqref{rel_stop}) then
\\\hspace*{2em}
stop iteration
\\\hspace*{1em}
endif
\\\hspace*{1em}
define $f_k(y)$, $A_k$, $\tilde{g}(t):=f_k(y_k(t))+g(t)$ 
\\\hspace*{1em}
set inner iteration tolerance $\texttt{tol}_{\mathtt{lin}}$
(cf.~\eqref{lin_stop},\eqref{lin_s2})
\\\hspace*{1em}
compute the low rank approximation~\eqref{Upt}, $\tilde{g}(t) \approx U p(t)$
\\\hspace*{1em}
solve linear IVP~\eqref{nonlRa} by inner block Krylov subspace iteration
(linear EBK)
\\\hspace*{1em}
$\Rightarrow$ obtain $y_{k+1}(t)=V_{\ell\cdot m} u_{k+1}(t)$, $t\in[0,T]$
\\
endfor (main iteration loop)
\hrule

\section{Numerical experiments}

\subsection{1D Burgers equation}
\label{s:num_exp1D}
This IVP is taken from~\cite[Ch.~IV.10]{ODE2(HW)}.
To find unknown function $u(x,t)$, we solve a 1D Burgers equation 
$$
u_t = \nu u_{xx} -uu_x , \quad 0\leqs x\leqs 1, 
\quad 0\leqs t\leqs T, 
$$
with viscosity $\nu = 3\cdot 10^{-4}$ and  $\nu = 3\cdot 10^{-5}$.
Initial and boundary conditions are
$$
u(x,0) = \frac32 x(1-x)^2, \quad u(0,t)=u(1,t)=0,
$$
and the final time $T$ is set to $0.5$, $1.0$, and $1.5$.
We discretize the Burgers equation by finite differences on a
uniform mesh with $N$ internal nodes $x_j=j\Delta x$, $j=1,\dots,n$,
$\Delta x=1/(N+1)$.  
To avoid numerical energy dissipation and to guarantee that a discretization of 
the advection term $uu_x$
results in  a skew-symmetric matrix~\cite{Krukier79}, we apply the standard second-order
central differences in the following way
$$
uu_x = \frac13 uu_x + \frac23 \left(\frac{u^2}2\right)_x\approx
\frac13 u_i\frac{u_{i+1}-u_{i-1}}{2\Delta x} + \frac23
\frac{u_{i+1}^2-u_{i-1}^2}{4\Delta x}.
$$
The diffusion term $\nu u_{xx}$ is discretized by second-order central differences
in the usual way.  Then this finite difference discretization
leads to an IVP of the form
\begin{equation}
\label{1D_discr}
y'(t) = -\Asymm y(t) - \Askew (y(t))y(t), \quad y(0)=v,  
\end{equation}
where 
$$
\Asymm y\approx -\nu u_{xx}, \quad \Askew (y)y \approx uu_x,
$$ 
$\Asymm =\Asymm^T\in\Rr^{N\times N}$ is positive definite and $(\Askew (y))^T=-\Askew (y)\in\Rr^{N\times N}$ for all $y\in\Rr^N$.
\added{It is important to make sure that discretization of the convective
term results in a skew-symmetric matrix~\cite{Krukier79,VerstappenVeldman2003}.
Otherwise, there are contributions of the convective term to
the symmetric part of the Jacobian matrix which distort the proper
energy dissipation in the ODE system and may prevent the 
the symmetric part from being positive semidefinite.
In the context of this paper, this would mean that $\omega\geqs 0$
is not guaranteed anymore.}

A straightforward way to apply the nonlinear iteration~\eqref{nonlR} to~\eqref{1D_discr}
would be setting in~\eqref{split_f} $A_k:=\Asymm $, $f_k(y):=-\Askew (y)y$.
However, to decrease the Lipschitz constant of $f_k(y)$ (cf.\ condition~\eqref{conv}), we also
include a linear part of $-\Askew (y)y$ to the $-A_ky$ term in~\eqref{split_f}.
More specifically, we apply the nonlinear iteration~\eqref{nonlR} to semidiscrete
Burgers equation~\eqref{1D_discr} in the form
\begin{equation}
\label{nonlR_1D}  
\begin{aligned}
y_{k+1}'(t) =
- (\,\underbrace{\Asymm  + \Askew (\bar{y}_k)}_{\displaystyle A_k}\,)y_{k+1}(t)
+ \underbrace{\left[\Askew (\bar{y}_k) - \Askew (y_k(t))\right]y_k(t)}_{\displaystyle f_k(y_k(t))},
\\
k=0,1,2,\dots,
\end{aligned}
\end{equation}
where $\bar{y}_k=y_k(T)$.
The nonlinear residual~\eqref{rk} in this case reads, for
$k=0,1,\dots$,
\begin{align*}
r_{k+1}(t) &= f_k(y_{k+1}(t)) - f_k(y_{k}(t))\\
          &= 
\left[\Askew (\bar{y}_{k}) - \Askew (y_{k+1}(t))\right]y_{k+1}(t) -
\left[\Askew (\bar{y}_{k}) - \Askew (y_k(t))\right]y_k(t),
\\
r_{k+1}(T) &= \left[\Askew (\bar{y}_k) - \Askew (y_{k+1}(T))\right]y_{k+1}(T).
\end{align*}
\added{Since the symmetric part of the matrix $A_k$ is 
$\frac12(A_k+A_k^T) = \Asymm$, the $\omega$ value is the smallest
eigenvalue of $\Asymm$.  The lowest eigenmodes should be well
approximated even on the coarsest grid.  Indeed, computing
$\omega$ numerically (which is easy for this one-dimensional
problem) yields approximately the same result on all grids:
we get $\omega\approx\texttt{2.96e-03}$ for $\nu=3\cdot 10^{-4}$
and $\omega\approx\texttt{2.96e-04}$ for $\nu=3\cdot 10^{-5}$.}
We stop our nonlinear iterative process~\eqref{nonlR} as soon as
\begin{equation}
\label{stop}
\|r_k(T)\|\leqslant \texttt{tol},
\end{equation} 
where $\texttt{tol}$ is a tolerance value.
As we primarily aim at solving IVPs stemming from space discretization of PDEs,
where space discretization error is significant,
we settle for a moderate tolerance value $\texttt{tol}=10^{-3}$
in this test.
To obtain the low rank representation~\eqref{Upt} at each iteration~$k$
we use $n_s=100$ sample points $t_j$, $j=1,\dots,n_s$ and set the number $m$ of 
singular vectors to $m=7$.
The largest truncated singular value $\sigma_{m+1}$, serving as an indicator
of the relative low rank representation error in~\eqref{Upt}, is then of order $10^{-6}$ 
for $T=0.5$ and of order $10^{-3}$ for $T=1.5$.

For the inner iterative solver EBK employed to solve~\eqref{nonlR} we set
the residual tolerance to $\texttt{tol}_{\mathtt{lin}}=\texttt{tol}$,
\added{which means that inner iterations are stopped as soon as
the inner linear residual $\rlin_{k+1}(t)$ satisfies
\begin{equation}
\label{lin_stop}
\|\rlin_{k+1}(T)\| \leqs \texttt{tol}_{\mathtt{lin}}.
\end{equation}}
Krylov subspace dimension is set to~$10$.   As the block size
in the Krylov subspace process is $m=7$, the EBK solver requires $m(10+1)=77$
vectors to store.
Since the problem becomes stiffer as the spatial grid gets finer, we use the EBK solver
in the shift-and-invert (SAI) mode~\cite{MoretNovati04,EshofHochbruck06}, i.e., 
the block Krylov subspace is computed
for the SAI matrix $(I+\gamma A_k)^{-1}$, $\gamma=T/10$, with $A_k$ defined in~\eqref{nonlR_1D}.
A sparse (banded in this test) LU~factorization of $I+\gamma A_k$ is computed once
at each nonlinear iteration and used every time an action of $(I+\gamma A_k)^{-1}$ is required.

To compare our nonlinear iterative EBK solver to another solver,
we also solve this test with a MATLAB stiff ODE solver \texttt{ode15s}.
This is a variable step size and variable order implicit multistep 
method~\cite{ODEsuite}.  The \texttt{ode15s} solver is run with
absolute and relative tolerances set, respectively, to
$\texttt{AbsTol} = \texttt{tol}$ and $\texttt{RelTol} = 10^4 \texttt{tol}$,
with $\texttt{tol}=10^{-9}$.
For these tolerance values both \texttt{ode15s} and our nonlinear EBK
deliver a comparable accuracy.  The relative error values reported below are computed 
as 
\begin{equation}
\label{err_reached}
\frac{\|y-y_{\mathrm{ref}}\|}{\|y_{\mathrm{ref}}\|},
\end{equation}
where $y$ is a numerical solution at final time $T$ and $y_{\mathrm{ref}}$ 
is a reference solution computed at final time by the \texttt{ode15s}
solver run with strict tolerance values 
$\texttt{AbsTol} = 10^{-12}$ and $\texttt{RelTol} = 10^{-8}$.
As $y_{\mathrm{ref}}$ is computed on the same spatial grid as $y$, the error value can be seen
as a reliable measure of time accuracy~\cite{RKC97}.

Actual delivered accuracies and corresponding 
required computational work for both solvers are reported in 
Table~\ref{t:Brg} (for viscosity $\nu = 3\cdot 10^{-4}$) 
and Table~\ref{t:Brg_mu3e_5} (for $\nu = 3\cdot 10^{-5}$).
For our nonlinear EBK solver the work is reported as the number of 
nonlinear iterations, with one LU~factorization and one matrix-vector
product (matvec) per iteration, and the total number of linear systems
solutions (which equals the total number of Krylov steps times 
the block size $m$).
The reported total number of 
LU~applications is not necessarily a multiple of the block size $m=7$ because
at the first nonlinear iteration the approximate solution is constant in time
and, hence, we set $m=1$ in~\eqref{Upt} at the first iteration.
The computational work for the \texttt{ode15s} solver is reported as the number 
of time steps, computed LU~factorizations and the ODE right hand side function 
evaluations (fevals).
\added{In Tables~\ref{t:Brg} and~\ref{t:Brg_mu3e_5} we also
report the norm of $\Asymm + \Askew(y(T))$ which can be seen
as a measure of the ODE system stiffness.}

\begin{table}
\caption{The 1D Burgers test problem, viscosity $\nu = 3\cdot 10^{-4}$.  
Attained error and 
computational work for our nonlinear EBK method and 
the \texttt{ode15s} solver. 
For the EBK method the work is reported as the nonlinear iteration number,
number of LU~factorizations, their applications and matvecs.
For the \texttt{ode15s} solver the work is measured the 
as the number of time steps, LU~factorizations, their applications and fevals.}
\label{t:Brg}  

\renewcommand{\arraystretch}{1.1} 
\begin{tabular}{ccccc}
\hline\hline
      &  method,            & iter./& LUs (LUs applic.), & relative\\ 
$T$   & tolerance           & steps & matvecs/fevals     & error \\
\hline
\multicolumn{5}{c}{grid $N=500$, $\|\Asymm  + \Askew (y(T))\|_1\approx 300$}\\
$0.5$  & nonlin.EBK($m=7$), {\tt 1e-03} & 5     & 5 (141), 5         & {\tt 5.17e-06} \\
$1.0$  & nonlin.EBK($m=7$), {\tt 1e-03} & 7     & 7 (220), 7         & {\tt 2.03e-05} \\
$1.5$  & nonlin.EBK($m=7$), {\tt 1e-03} & 10    & 10 (340), 10       & {\tt 5.31e-05} \\
$0.5$  & {\tt ode15s},  {\tt 1e-09}     & 59    & 14 (73), 575       & {\tt 1.93e-06} \\
$1.0$  & {\tt ode15s},  {\tt 1e-09}     & 70    & 16 (86), 588       & {\tt 8.38e-06} \\
$1.5$  & {\tt ode15s},  {\tt 1e-09}     & 83    & 19 (118), 620      & {\tt 2.21e-05} \\
\hline
\multicolumn{5}{c}{grid $N=1000$, $\|\Asymm  + \Askew (y(T))\|_1\approx 1\,200$}\\
$0.5$  & nonlin.EBK($m=7$), {\tt 1e-03} & 5     & 5 (170), 5         & {\tt 5.06e-06} \\
$1.0$  & nonlin.EBK($m=7$), {\tt 1e-03} & 7     & 7 (256), 7         & {\tt 2.00e-05} \\
$1.5$  & nonlin.EBK($m=7$), {\tt 1e-03} & 10    & 10 (389), 10       & {\tt 5.30e-05} \\
$0.5$  & {\tt ode15s},  {\tt 1e-09}     & 67    & 15 (83), 1085      & {\tt 1.76e-06} \\
$1.0$  & {\tt ode15s},  {\tt 1e-09}     & 78    & 17 (104), 1106     & {\tt 1.13e-05} \\
$1.5$  & {\tt ode15s},  {\tt 1e-09}     & 91    & 19 (134), 1136     & {\tt 2.48e-05} \\
\hline
\multicolumn{5}{c}{grid $N=2000$, $\|\Asymm  + \Askew (y(T))\|_1\approx 4\,800$}\\
$0.5$  & nonlin.EBK($m=7$), {\tt 1e-03} & 5     & 5 (177), 5         & {\tt 5.07e-06} \\
$1.0$  & nonlin.EBK($m=7$), {\tt 1e-03} & 7     & 7 (277), 7         & {\tt 2.00e-05} \\
$1.5$  & nonlin.EBK($m=7$), {\tt 1e-03} & 11    & 11 (452), 11       & {\tt 4.38e-05} \\
$0.5$  & {\tt ode15s},  {\tt 1e-09}     & 73    & 16 (91), 2093      & {\tt 2.29e-06} \\
$1.0$  & {\tt ode15s},  {\tt 1e-09}     & 85    & 19 (107), 2109     & {\tt 7.96e-06} \\
$1.5$  & {\tt ode15s},  {\tt 1e-09}     & 98    & 22 (139), 2141     & {\tt 2.22e-05} \\
\hline
\multicolumn{5}{c}{grid $N=4000$, $\|\Asymm  + \Askew (y(T))\|_1\approx 19\,200$}\\
$0.5$  & nonlin.EBK($m=7$), {\tt 1e-03} & 5    & 5 (193), 5         & {\tt 5.06e-06} \\
$1.0$  & nonlin.EBK($m=7$), {\tt 1e-03} & 8    & 8 (347), 8         & {\tt 4.82e-06} \\
$1.5$  & nonlin.EBK($m=7$), {\tt 1e-03} & 11   & 11 (501), 11       & {\tt 4.38e-05} \\
$0.5$  & {\tt ode15s},  {\tt 1e-09}     & 79   & 18 (99), 4101      & {\tt 1.93e-06} \\
$1.0$  & {\tt ode15s},  {\tt 1e-09}     & 90   & 20 (112), 4114     & {\tt 8.25e-06} \\
$1.5$  & {\tt ode15s},  {\tt 1e-09}     & 103  & 23 (144), 4146     & {\tt 2.21e-05} \\
\hline 
\end{tabular}
\end{table}

\begin{figure}
\includegraphics[width=0.49\textwidth]{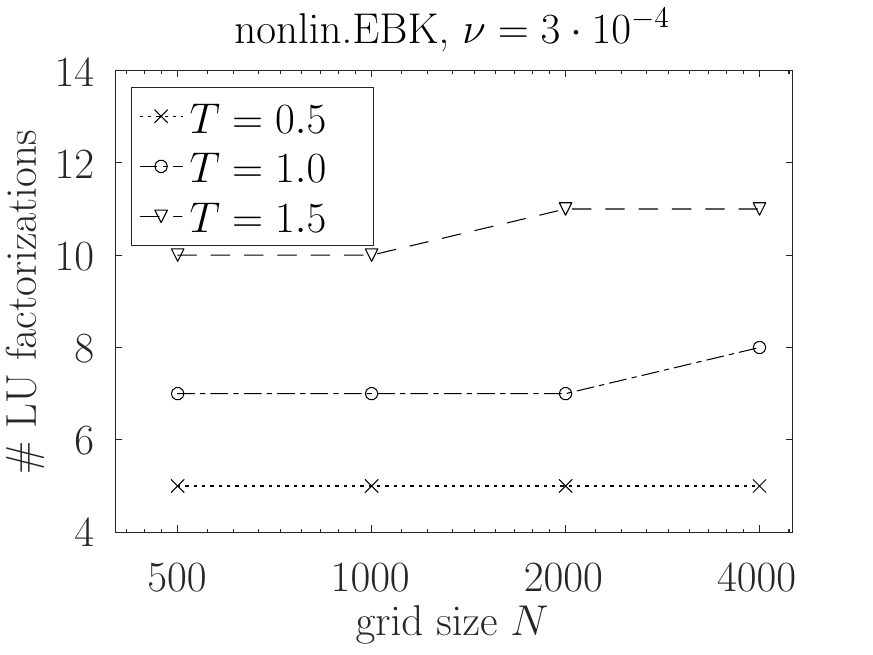}%
\includegraphics[width=0.49\textwidth]{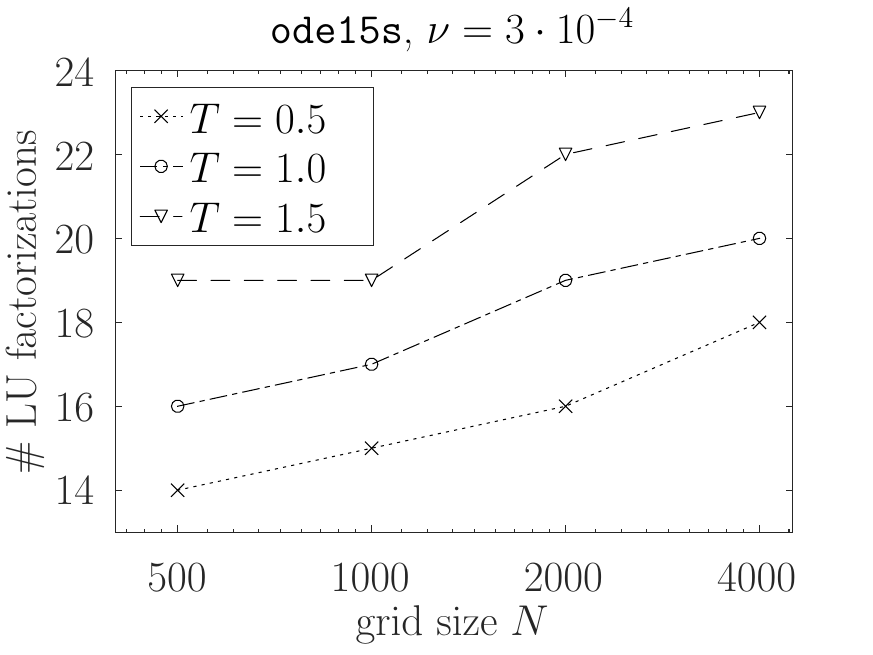}  
\caption{The Burgers test.
Number of required LU~factorizations versus the grid size 
in the nonlinear EBK solver (left) and \texttt{ode15s} (right) for
different final time values $T$ and viscosity $\nu=3\cdot 10^{-4}$}
\label{f:LUs}  
\end{figure}

As we can see in Tables~\ref{t:Brg} and~\ref{t:Brg_mu3e_5}, 
our solver requires less LU~factorizations than the \texttt{ode15s} solver.
It does require more LU~factorization actions but these are relatively cheap.
Moreover, in EBK they are carried out simultaneously for blocks of $m$ right hand
sides.
For both solvers, in Figure~\ref{f:LUs} we plot the numbers of LU factorizations
(presented in the tables) versus the grid size.
The number of nonlinear iterations in the EBK solver
(with one LU~factorization per iteration) remains practically constant as
the grid gets finer, whereas the number of LU~factorizations in \texttt{ode15s}
increases with the grid size.  
Furthermore, we see that the \texttt{ode15s} solver requires
more time steps and significantly more fevals on the finer grids.
$T=1.5$ is approximately the largest possible time for which our nonlinear EBK 
solver converges in this test.  From Figure~\ref{f:LUs} we see that 
a higher efficiency is reached for larger $T$ values.
Indeed, on the finest grid $N=4000$ 
$5\times 2= 10$~LU factorizations per unit time are required for $T=0.5$,
$8$~LU factorizations for $T=1$ and $11/1.5\approx 7$ factorizations for
$T=1.5$.    

Comparing Tables~\ref{t:Brg} and~\ref{t:Brg_mu3e_5}, we see that performance of 
the \texttt{ode15s} solver improves for the smaller viscosity value $\nu=3\cdot 10^{-5}$.
This is expected as the stiffness of the space-discretized Burgers equation
decreases with $\nu$.  Indeed, on finer grids, where contributions of $\Asymm$ 
(proportional to $(\Delta x)^{-2}$) dominate those of $\Askew(t)$ 
(proportional to $(\Delta x)^{-1}$),
we have $\|\Asymm+\Askew(t)\|\approx\mathcal{O}(\nu)$, see the values
$\|\Asymm+\Askew(T)\|$ reported in the tables. 

\begin{figure}  
\centerline{\includegraphics[width=0.6\textwidth]{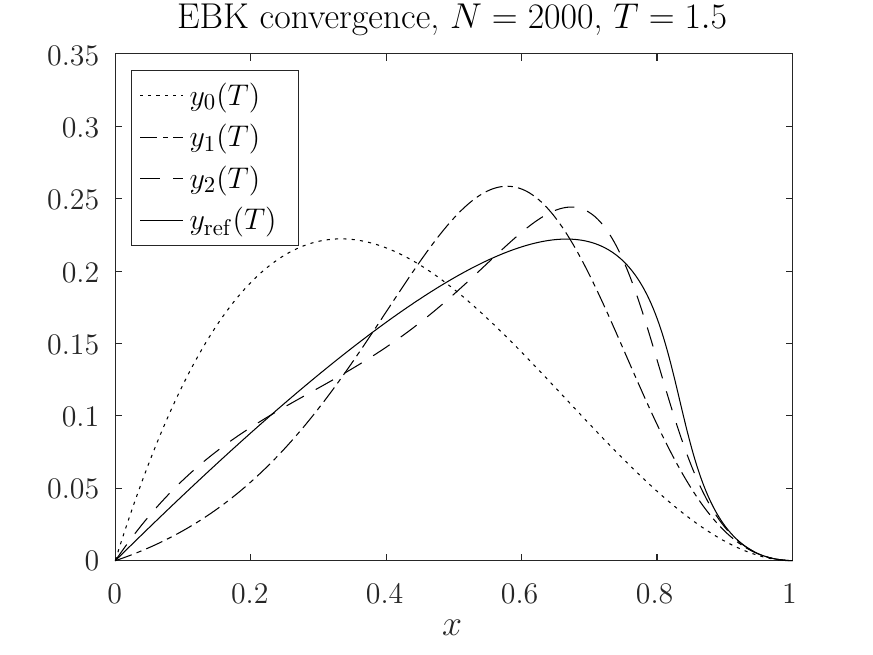}}
\caption{The Burgers test.
The first three iterative approximations $y_0(T)=v$, $y_1(T)$, $y_2(T)$
and reference solution $y_{\mathrm{ref}}(T)$\deleted{ for the Burgers test}, 
grid size $n=2000$, final time $T=1.5$, viscosity $\nu=3\cdot 10^{-4}$.}
\label{f:1D_Brgrs}
\end{figure}

\begin{figure}
\includegraphics[width=0.48\textwidth]{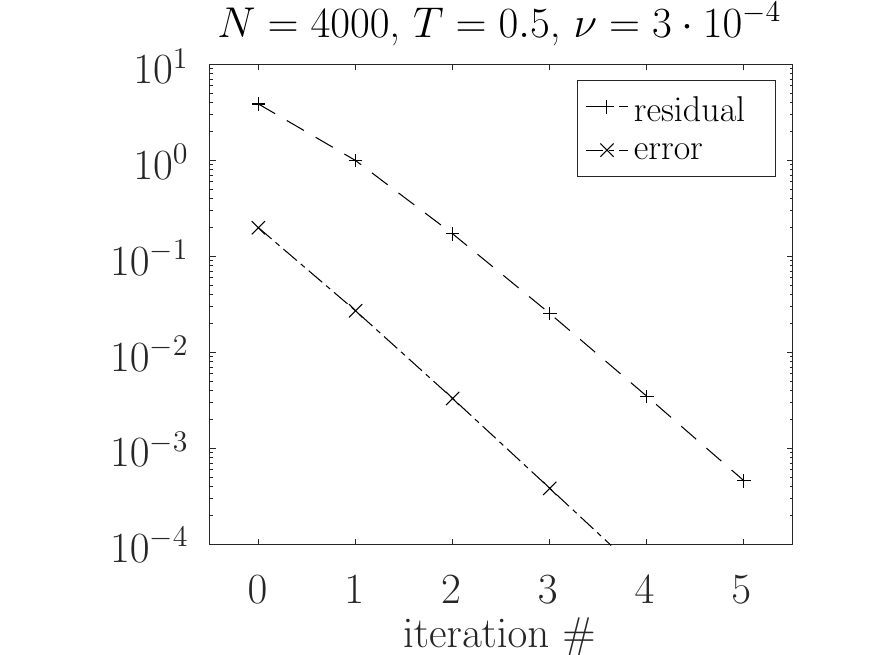}%
\includegraphics[width=0.48\textwidth]{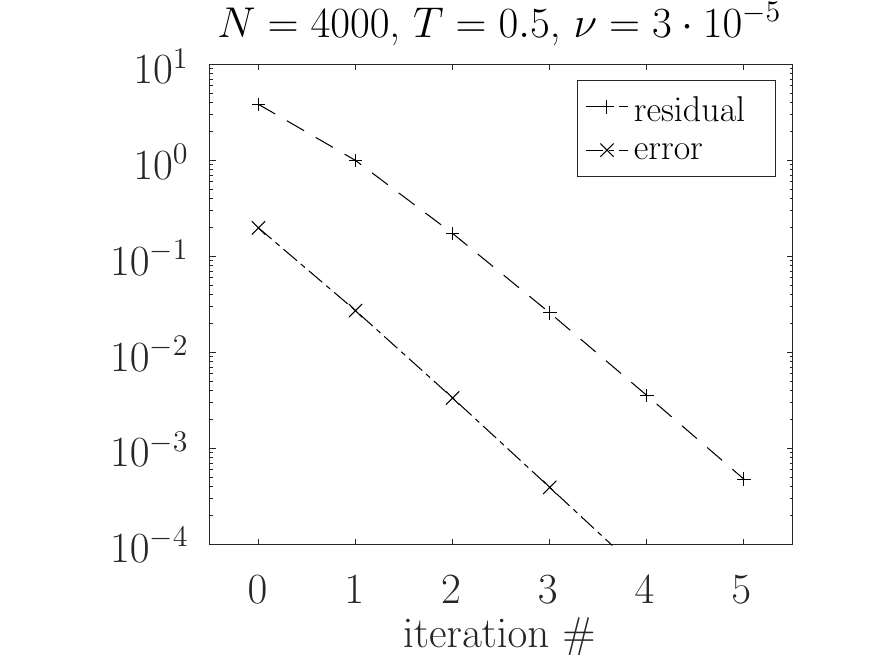}\\[2ex]
\includegraphics[width=0.48\textwidth]{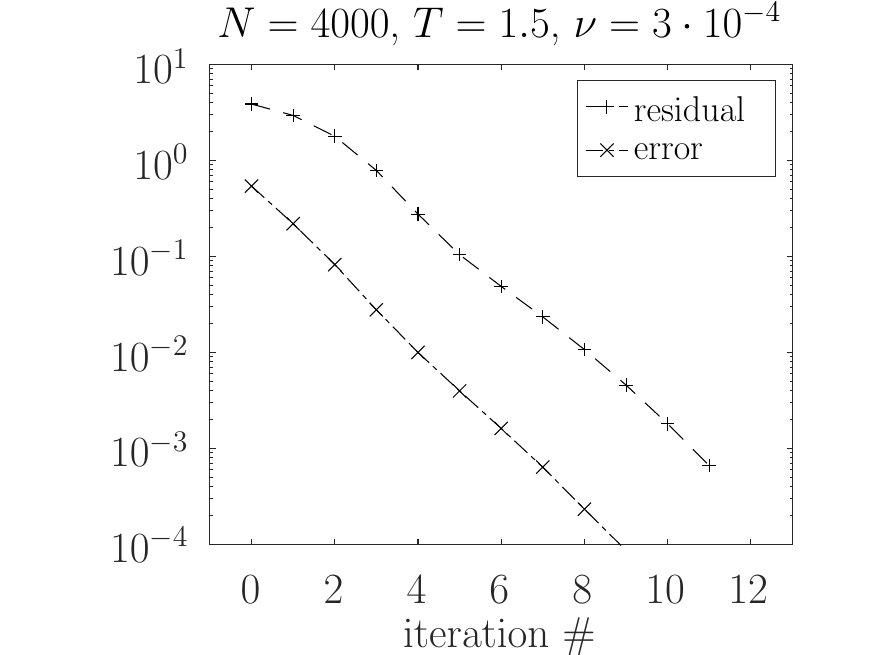}%
\includegraphics[width=0.48\textwidth]{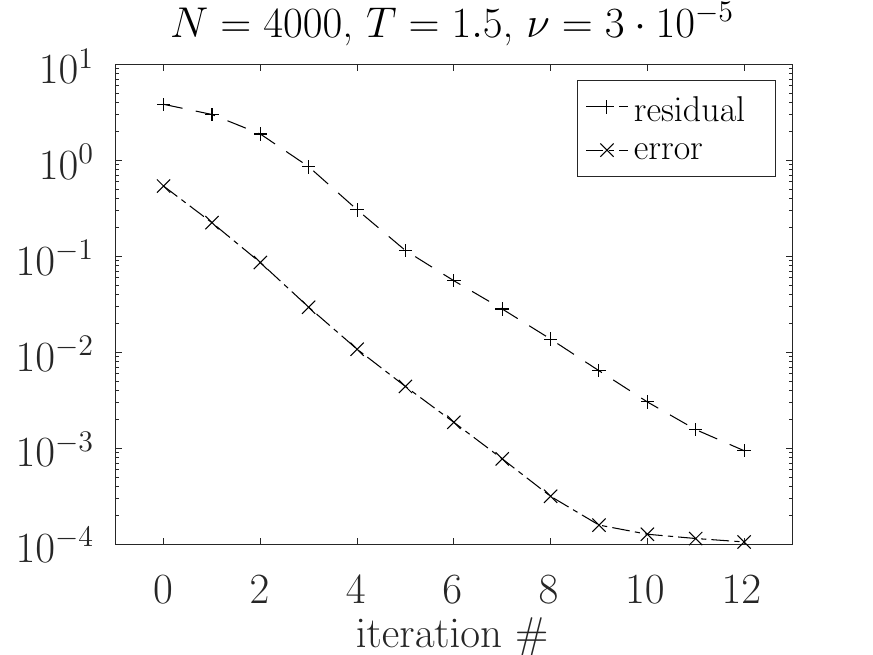}
\caption{The Burgers test.
Residual norm~\eqref{rk} and relative error norm~\eqref{err_reached} versus 
iteration number for $T=0.5$ (top plots) and $T=1.5$ (bottom plots),
$\nu=3\cdot 10^{-4}$ (left plots) and $\nu=3\cdot 10^{-5}$ (right plots).
The grid size is $n=4000$. Error stagnation visible at the bottom right plot
can be repaired by slightly increasing the block size $m$.}
\label{f:conv}  
\end{figure}

In Figure~\ref{f:1D_Brgrs} we plot the reference solution $y_{\mathrm{ref}}(T)$ and 
the first three iterative approximations $y_0(T)=v$, $y_1(T)$ and $y_2(T)$.
Convergence plots of the residual and error norms are shown in Figure~\ref{f:conv}.
There, the shown residual norm is computed at $t=T$ according to~\eqref{rk} and
the error norm is the relative norm defined in~\eqref{err_reached}.
As we see, the residual turns out to be a reliable measure of the error.
Furthermore, convergence behavior is not influenced by viscosity~$\nu$, which is
in accordance with theory:
$\nu$ does not change the nonlinear term $f_k$ in~\eqref{nonlR}
and its Lipschitz constant.

\begin{table}
\caption{The 1D Burgers test problem, viscosity $\nu = 3\cdot 10^{-5}$.  
Attained error and computational work for our nonlinear EBK method and 
the \texttt{ode15s} solver. 
For the EBK method the work is reported as the nonlinear iteration number,
number of LU~factorizations, their applications and matvecs.
For the \texttt{ode15s} solver the work is measured the 
as the number of time steps, LU~factorizations, their applications and fevals.}
\label{t:Brg_mu3e_5}  

%
%
\renewcommand{\arraystretch}{1.1} 
\begin{tabular}{ccccc}
\hline\hline
      &  method,            & iter./& LUs (LUs applic.), & relative\\ 
$T$   & tolerance           & steps & matvecs/fevals     & error \\
\hline
\multicolumn{5}{c}{grid $N=500$, $\|\Asymm  + \Askew (y(T))\|_1\approx 90$}\\
$0.5$  & nonlin.EBK($m=7$), {\tt 1e-03} & 5     & 5 (69),  5         & {\tt 1.82e-05} \\
$1.0$  & nonlin.EBK($m=7$), {\tt 1e-03} & 7     & 7 (139), 7         & {\tt 2.26e-05} \\
$1.5$  & nonlin.EBK($m=7$), {\tt 1e-03} & 13    & 13 (414), 13       & {\tt 1.10e-04} \\
$0.5$  & {\tt ode15s},  {\tt 1e-09}     & 28    & 7 (38), 540        & {\tt 3.71e-06} \\
$1.0$  & {\tt ode15s},  {\tt 1e-09}     & 38    & 9 (53), 555        & {\tt 1.23e-05} \\
$1.5$  & {\tt ode15s},  {\tt 1e-09}     & 52    & 13 (89), 591       & {\tt 3.38e-05} \\
\hline
\multicolumn{5}{c}{grid $N=1000$, $\|\Asymm  + \Askew (y(T))\|_1\approx 210$}\\
$0.5$  & nonlin.EBK($m=7$), {\tt 1e-03} & 5     & 5 (90), 5          & {\tt 6.20e-06} \\
$1.0$  & nonlin.EBK($m=7$), {\tt 1e-03} & 7     & 7 (176), 7         & {\tt 2.25e-05} \\
$1.5$  & nonlin.EBK($m=7$), {\tt 1e-03} & 12    & 12 (430), 12       & {\tt 1.07e-04} \\
$0.5$  & {\tt ode15s},  {\tt 1e-09}     & 35    & 9  (44), 1046      & {\tt 2.44e-06} \\
$1.0$  & {\tt ode15s},  {\tt 1e-09}     & 45    & 10 (60), 1062      & {\tt 1.64e-05} \\
$1.5$  & {\tt ode15s},  {\tt 1e-09}     & 59    & 15 (99), 2102      & {\tt 3.83e-05} \\
\hline
\multicolumn{5}{c}{grid $N=2000$, $\|\Asymm + \Askew (y(T))\|_1\approx 540$}\\
$0.5$  & nonlin.EBK($m=7$), {\tt 1e-03} & 5     & 5 (120), 5         & {\tt 5.29e-06} \\
$1.0$  & nonlin.EBK($m=7$), {\tt 1e-03} & 7     & 7 (190), 7         & {\tt 2.22e-05} \\
$1.5$  & nonlin.EBK($m=7$), {\tt 1e-03} & 12    & 12 (494), 12       & {\tt 1.06e-04} \\
$0.5$  & {\tt ode15s},  {\tt 1e-09}     & 42    & 10 (55), 2057      & {\tt 3.91e-06} \\
$1.0$  & {\tt ode15s},  {\tt 1e-09}     & 52    & 12 (70), 2072      & {\tt 1.32e-05} \\
$1.5$  & {\tt ode15s},  {\tt 1e-09}     & 66    & 17 (108), 4111     & {\tt 3.32e-05} \\
\hline
\multicolumn{5}{c}{grid $N=4000$, $\|A_0 + \Askew (y(T))\|_1\approx 1\,920$}\\
$0.5$  & nonlin.EBK($m=7$), {\tt 1e-03} & 5    & 5 (149), 5         & {\tt 5.24e-06} \\
$1.0$  & nonlin.EBK($m=7$), {\tt 1e-03} & 8    & 8 (276), 8         & {\tt 5.52e-06} \\
$1.5$  & nonlin.EBK($m=7$), {\tt 1e-03} & 12   & 12 (578), 12       & {\tt 1.07e-04} \\
$0.5$  & {\tt ode15s},  {\tt 1e-09}     & 47   & 11 (61), 4063      & {\tt 3.86e-06} \\
$1.0$  & {\tt ode15s},  {\tt 1e-09}     & 58   & 15 (80), 4082      & {\tt 1.22e-05} \\
$1.5$  & {\tt ode15s},  {\tt 1e-09}     & 71   & 18 (112), 8115     & {\tt 3.36e-05} \\
\hline 
\end{tabular}
\end{table}

\subsection{3D Bratu test problem}
\label{s:num_exp3D}
\begin{figure}
\centerline{%
\includegraphics[width=0.47\linewidth]{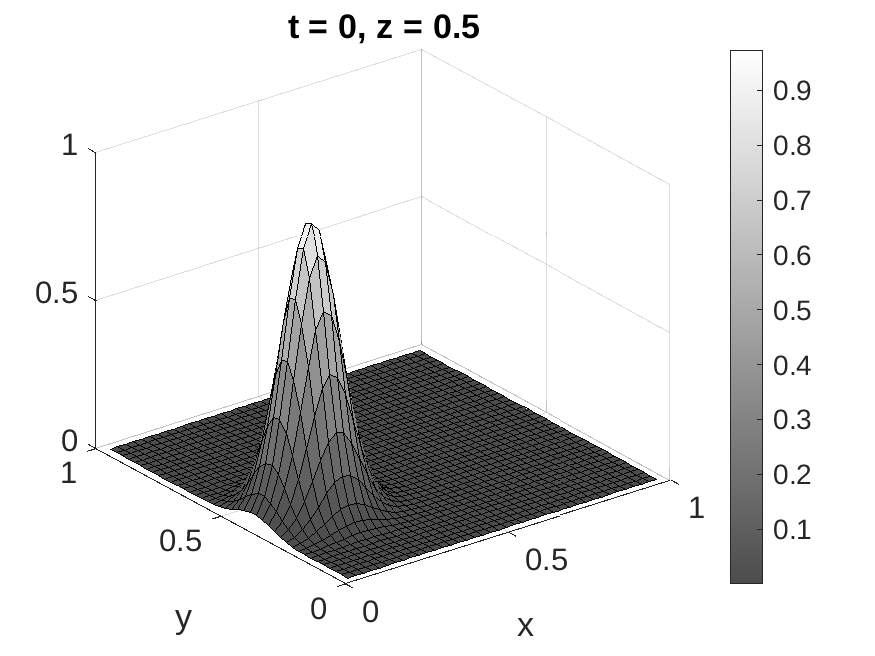}\hfill
\includegraphics[width=0.47\linewidth]{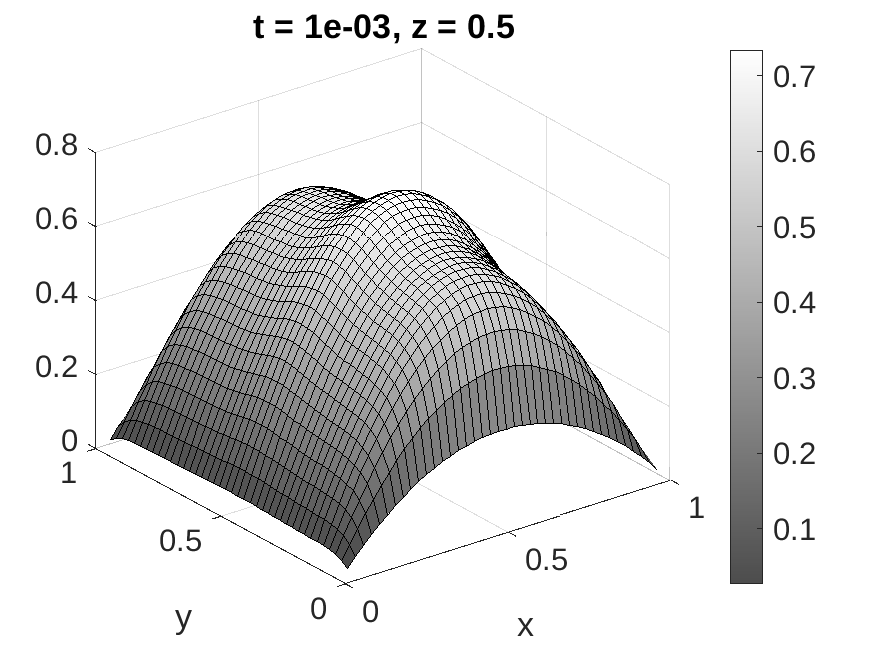}}  
\caption{Left: initial value function $u(x,y,z,0)$ of the Bratu test 
on a uniform $40\times 40\times 40$ grid for $z=0.5$.  
Right: numerical solution $u(x,y,z,t)$ on the same grid for 
$z=0.5$ and $t=1\cdot 10^{-3}$.}
\label{f:3D_Bratu}  
\end{figure}

We now consider a nonstationary IVP for Liouville--Bratu--Gelfand equation~\cite{JacobsenSchmitt2002}:
find $u=u(x,y,z,t)$ such that
\begin{equation}
\label{Bratu}
\begin{aligned}
& u_t = 10^4u_{xx} + 10^2u_{yy} + u_{zz} + Ce^u + g^{\mathsf{src}}(x,y,z,t), 
\\ 
& \qquad\qquad (x,y,z)\in\Omega=[0,1]^3,  \quad 0\leqs t\leqs T, \quad C=3\cdot 10^4,
\\
& u(x,y,z,0) = e^{-100\left( (x-0.2)^2 + (y-0.4)^2 + (z-0.5)^2 \right)},\quad
u\bigl|_{\partial\Omega}=0.
\end{aligned}
\end{equation}
\added{Here} the source function $g^{\mathsf{src}}$ is defined as
\begin{gather*}
g^{\mathsf{src}}(x,y,z,t) = 
 \begin{cases}
    \added{e^{-100((x-x_0(t))^2 + (y-y_0(t))^2 + (z-0.5)^2)}} + 
    C u(x,y,z,0), \; & t\leqslant 5\cdot 10^{-5},\\
    \added{e^{-100((x-x_0(t))^2 + (y-y_0(t))^2 + (z-0.5)^2)}},            
                  \; & t>         5\cdot 10^{-5},
  \end{cases}
\\
\added{x_0(t) = 0.5 + 0.3\cos(2000\pi t),\quad
y_0(t) = 0.5 + 0.3\sin(2000\pi t),}
\end{gather*}
and the final time $T$ is either $T=5\cdot 10^{-5}$ or $T=1\cdot 10^{-4}$.
Solution to~\eqref{Bratu} is shown in Figure~\ref{f:3D_Bratu}.

The regular second order finite-difference discretization of the spatial derivatives
in~\eqref{Bratu} on a uniform grid yields an IVP
\begin{equation}
\label{3D_discr}
y'(t) = -Ay(t) + \hat{f}(y(t)) + g(t), \quad y(0)=v,  
\end{equation}
where the entries of the vector functions $y(t)$ and $g(t)$ contain the values
of the numerical solution and function $g^{\mathsf{src}}(x,y,z,t)$ at the grid nodes,
respectively, and
\begin{equation}
\label{anis}  
A y(t)\approx -(10^4u_{xx} + 10^2u_{yy} + u_{zz}).
\end{equation}
Furthermore, the components $[\hat{f}(y)]_i$ of the vector function $\hat{f}(y)$ are
defined as 
$$
\left[\hat{f}(y)\right]_i = C e^{y_i},\quad i=1,\dots,N.
$$
To apply nonlinear waveform relaxation~\eqref{nonlR} for solving~\eqref{3D_discr},
one could have chosen to set $A_k:=A$ and $f_k(y):=\hat{f}(y)$. 
However, to have a smaller Lipschitz constant $L$ in $f_k(y)$, we supplement
$A_k$ with a linear part of $\hat{f}$, setting
\begin{equation}
\label{nonlR_3D}
\begin{aligned}
y_{k+1}'(t) =
- (\,\underbrace{A - J(\bar{y}_k)}_{\displaystyle A_k}\,)y_{k+1}(t)
+ \underbrace{\left[\hat{f}(y_k(t)) - J(\bar{y}_k)y_k(t)\right]}_{\displaystyle f_k(y_k(t))} + g(t),
\\
k=0,1,2,\dots,
\end{aligned}
\end{equation}
where $J(\bar{y}_k)$ is the Jacobian matrix of $\hat{f}(y)$ evaluated at 
$\bar{y}_k=y_k(T)$ and an approximation
$$
\hat{f}(y_{k+1}(t)) \approx \hat{f}(y_k(t)) + J(\bar{y}_k) (y_{k+1}(t)-y_k(t))
$$
is used.  The nonlinear residual~\eqref{rk} takes in this case the form, 
for $k=0,1,\dots$,
\begin{align*}
r_{k+1}(t) &= f_k(y_{k+1}(t)) - f_k(y_{k}(t))\\
          &= 
\left[\hat{f}(y_{k+1}(t)) - J(\bar{y}_k)y_{k+1}(t)\right] -
\left[\hat{f}(y_k(t)) - J(\bar{y}_k)y_k(t)\right]
\\
&= \hat{f}(y_{k+1}(t)) - \hat{f}(y_k(t)) - J(\bar{y}_k)(y_{k+1}(t) - y_k(t)).
\end{align*}
\deleted{The iterations are stopped provided that
$\|r_k(T)\|\leqs\texttt{tol}$ for the tolerance \texttt{tol} set in this test to
$10^{-2}$.}
\added{Since this IVP is stiff and anisotropic, the nonlinear residual
norm can be very large in norm, $\mathcal{O}(10^8)$.  
Therefore we switch from the absolute stopping criterion 
$\|r_k(T)\|\leqs\texttt{tol}$ used in the previous test
to a relative stopping criterion
\begin{equation}
\label{rel_stop}
\|r_k(T)\|\leqs\texttt{tol} \|r_0(T)\|,
\end{equation}
with tolerances \texttt{tol} varying in this test from
$10^{-2}$ to $10^{-4}$ (the specific values are reported
in Table~\ref{t:3D_Bratu} discussed below).
Note that $y_0(t)\equiv v$, as discussed in Section~\ref{s:implem}.}
 
The low rank representation~\eqref{Upt} at each iteration~$k$ is computed
at $n_s=100$ sample points $t_j$, $j=1,\dots,n_s$ for either $m=4$ or $m=5$
low rank terms.
The \deleted{residual tolerance in the} 
EBK solver applied to solve~\eqref{nonlR} is 
\replaced{stopped as soon as}{defined as}
\begin{equation}
\label{lin_s2}
\|\rlin_{k+1}(T)\| \leqs 
\texttt{tol}_{\mathtt{lin}}=\|f_k(y_k(0)) + g(0)\| \, \texttt{tol} /10,  
\end{equation}
\added{which means that the linear residual should be small
with respect to the other terms in the perturbed ODE~\eqref{pertODE}.
We switch from the absolute inner stopping criterion~\eqref{lin_stop}
to~\eqref{lin_s2} because in this stiff anisotropic problem 
the term $f_k(y_k(t)) + g(t)$ can be very large in norm.}
Krylov subspace dimension is set to~$10$.
Hence, for the block size $m=5$, the EBK solver has to store $5\cdot (10+1)=55$
vectors.
Since the problem is stiff, just as in the previous test, the EBK solver is employed in the
SAI mode, with a sparse LU~factorization of $I+\gamma A_k$ computed and applied
throughout each nonlinear iteration.
The relative error values we report below are computed according to~\eqref{err_reached}.

We compare our nonlinear EBK solver to the two-stage Rosenbrock method
ROS2, see~\cite{ros2} or~\cite[Chapter IV.5]{HundsdorferVerwer:book}.
For IVP~\eqref{IVP0}, it can be written as
\begin{equation}
\label{ros2}
\begin{aligned}
y^{\ell+1} &= y^\ell + \frac32 \tau k_1 + \frac12 \tau k_2, \\
(I-\gamma\tau\widehat{A})k_1 &= \Phi(t_\ell,y^\ell), \\
(I-\gamma\tau\widehat{A})k_2 &= \Phi(t_{\ell+1},y^\ell + \tau k_1) - 2k_1.
\end{aligned}
\end{equation}
Here $y^\ell\approx y(\ell\tau)$ is numerical solution at time step $\ell=0,1,2,\dots$
($y^0=v$)
obtained with a time step size $\tau>0$.
This scheme is second order consistent for any matrix $\widehat{A}\in\Rr^{N\times N}$.
Here we take $\widehat{A}$ to be the Jacobian matrix of $\Phi$ computed at $t_\ell=\ell\tau$
and $y_\ell$.  At each time step a sparse LU factorization of $I-\gamma\tau\widehat{A}$
can be computed and used to solve the two linear systems in $k_{1,2}$.
The numerical tests presented here are carried in Matlab and, as it turns out,
for the grid sizes used,
it is much faster to compute the action of the matrix $(I-\gamma\tau\widehat{A})^{-1}$
by the Matlab backslash operator~\verb|\| than to compute a sparse LU factorization
first and to apply it twice\footnote{To solve a linear system 
\texttt{Ax=b} with a square nonsingular matrix \texttt{A}, 
one can use the Matlab backslash operator~$\mathtt{\backslash}$ as 
\texttt{x=A$\mathtt{\backslash}$b},
in which case the operator computes and applies
an LU~factorization.}.  
This is apparently because Matlab creates 
an overhead to call sparse LU factorization routines 
(of the UMFPACK software~\cite{UMFPACK1,UMFPACK2}) and then to store the computed 
LU factors as Matlab variables. 
Note that both the backslash operator as well as the LU factorization for large
sparse matrices in Matlab are based on UMFPACK.

\begin{figure}
\includegraphics[width=0.48\linewidth]{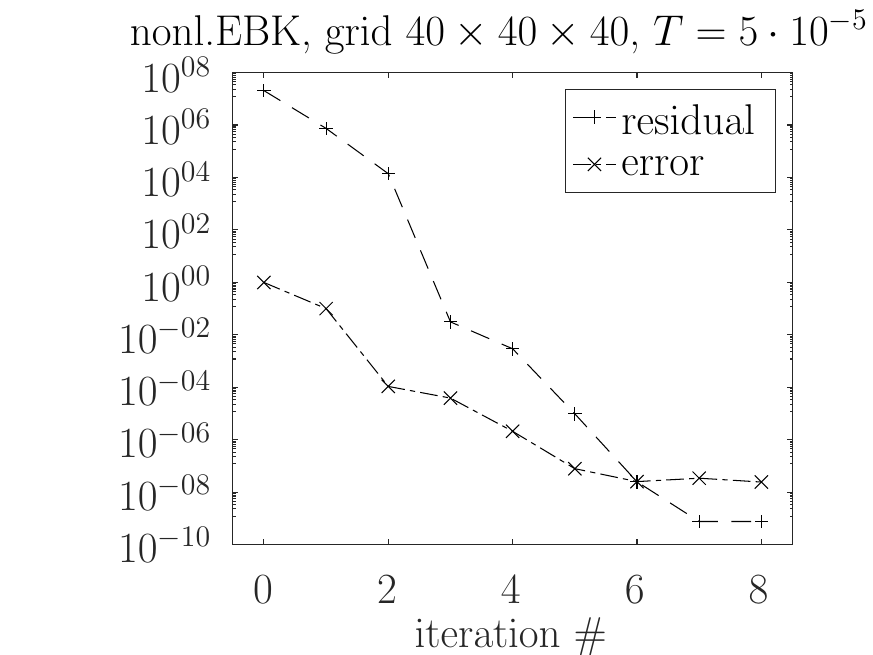}\hfill
\includegraphics[width=0.48\linewidth]{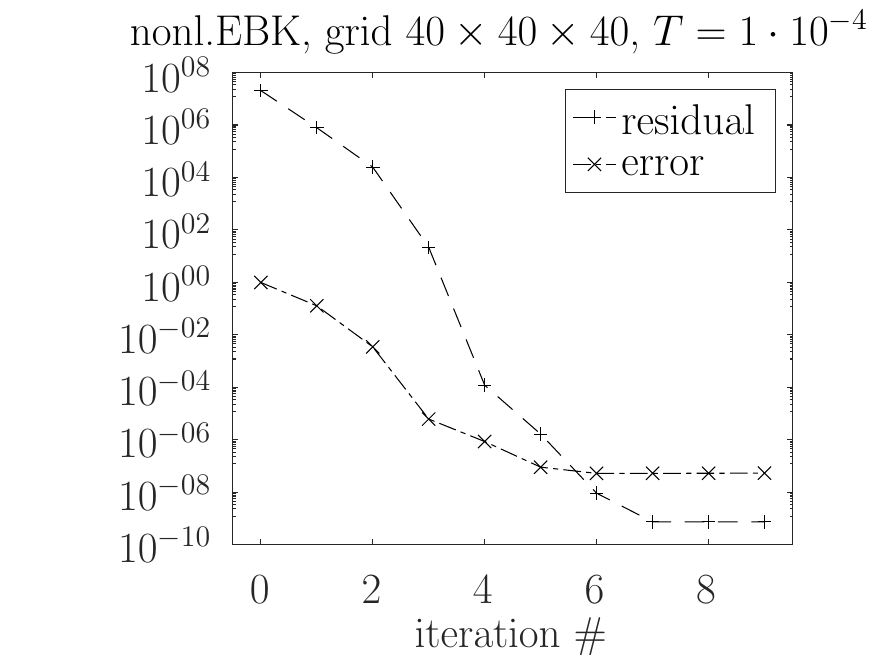}
\caption{The Bratu 3D test. 
Residual and \added{relative} error norm versus iteration number in the nonlinear 
EBK solver for the grid size $40\times 40\times 40$, 
\deleted{(top plots) and $60\times 60\times 60$ (bottom plots),}
$T=5\cdot10^{-5}$ (left plots) and
$T=1\cdot10^{-4}$ (right plots).
The plots are made with \added{increased block size $m=32$
and number of samples $n_s=300$.
Similar behavior is observed for the $60\times 60\times 60$ grid.}}
\label{f:ebk_Bratu}  
\end{figure}

Comparison results of our nonlinear EBK solver and ROS2 method are presented in
Table~\ref{t:3D_Bratu}.  The results clearly show that, for this test problem,
the nonlinear EBK solver is much more efficient than the ROS2 method both
in terms of required LU~factorizations (and their actions) and the CPU time.
We also see that the number of nonlinear waveform relaxation iterations 
\replaced{hardly}{does not} change\added{s}
with grid size and with final time $T$
\added{(see also Table~\ref{t:Bratu_conv})}.
While the convergence independence on the grid size is expected (as the 
grid size does not affect the nonlinear part $Ce^u$ in the ODE),
the weak convergence dependence on $T$ is probably a property
of the specific test problem.  From plots in Figure~\ref{f:ebk_Bratu}
we see that convergence behavior does \replaced{depend on}{change with} $T$. 

\begin{table}
\caption{The 3D Bratu test problem.  Attained error and 
computational work for our nonlinear EBK method and 
the ROS2 solver. 
For the EBK method the work is reported as the nonlinear iteration number,
number of LU~factorizations, their applications and matvecs.
For the ROS2 solver the work is shown as 
as the number of time steps, LU~factorizations, their applications and fevals.}
\label{t:3D_Bratu}  

\renewcommand{\arraystretch}{1.1} 
\begin{tabular}{ccccc}
\hline\hline
 method, time steps (t.s.) 
                     & iter./ & LUs (LUs applic.), & relative & CPU \\ 
 tolerance / $\tau$  & steps  & matvecs/fevals     & error    & time, s\\
\hline
\multicolumn{5}{c}{grid $40\times 40\times 40$, $T=5\cdot 10^{-5}$, $T\|A(y(T))\|_1\approx\texttt{3.40e+03}$}\\
nonlin.EBK($m=4$), 1 t.s., {\tt1e-02} &   2   &  2 (37), 2 & {\tt1.17e-04}  & 6.3 \\                
nonlin.EBK($m=5$), 1 t.s., {\tt1e-04} &   3   &  3 (111), 3 & {\tt4.04e-05}  & 12.9 \\
ROS2, $\tau = T/320$ & 320   & 320 (640), 640 & {\tt2.36e-04}  & 225 \\
ROS2, $\tau = T/640$ & 640   & 640 (1280), 1280 & {\tt5.73e-07}  & 442
%
%
%
\\\hline
\multicolumn{5}{c}{grid $40\times 40\times 40$, $T=1\cdot 10^{-4}$, $T\|A(y(T))\|_1\approx\texttt{6.79e+03}$}\\
%
nonlin.EBK($m=4$), 1 t.s., {\tt1e-03} & 3   & 3 (70), 3  & {\tt2.09e-05} & 10.4 \\
nonlin.EBK($m=5$), 1 t.s., {\tt1e-03} & 3   & 3 (80), 3  & {\tt1.38e-05} & 11.1 \\
ROS2, $\tau = T/320$              & 320 &  320 (640), 640 & {\tt1.51e-05}  & 201 \\
ROS2, $\tau = T/640$              & 320 & 640 (1280), 1280 & {\tt7.51e-06}  & 401
\\\hline
\multicolumn{5}{c}{\added{grid $40\times 40\times 40$, $T=1\cdot 10^{-3}$, $T\|A(y(T))\|_1\approx\texttt{6.79e+04}$}}\\
%
\added{nonlin.EBK($m=5$), 20 t.s., {\tt1e-03}} 
                                  & 46  & 46 (908), 46    & {\tt1.34e-05} & 148 \\
\added{nonlin.EBK($m=5$), 10 t.s., {\tt1e-03}} 
                                  & 30  & 30 (646), 30    & {\tt1.19e-05} & 102 \\
\added{ROS2, 3200 t.s., $\tau=T/3200$}
                                & 3200 & 3200 (6400), 3200 & {\tt2.19e-06} & 2007
\\\hline
\multicolumn{5}{c}{grid $60\times 60\times 60$, $T=5\cdot 10^{-5}$, $T\|A(y(T))\|_1\approx\texttt{7.50e+03}$}\\
%
nonlin.EBK($m=4$), 1 t.s., {\tt1e-02} & 2   & 2 (41), 2  & {\tt1.18e-04} & 51.5 \\
nonlin.EBK($m=4$), 1 t.s., {\tt1e-04} & 3   & 3 (91), 3  & {\tt6.91e-05} & 90 \\
%
ROS2, $\tau = T/320$              & 320 &  320 (640), 640 & {\tt2.35e-04}  & 2177 \\
%
ROS2, $\tau = T/640$              & 640 & 640 (1280), 1280 & {\tt1.66e-05}  & 4790
\\\hline
\multicolumn{5}{c}{grid $60\times 60\times 60$, $T=1\cdot 10^{-4}$, $T\|A(y(T))\|_1\approx\texttt{1.50e+04}$}\\
nonlin.EBK($m=4$), 1 t.s., {\tt1e-02} & 2   & 2 (38), 2  & {\tt3.48e-03} & 51\\
nonlin.EBK($m=4$), 1 t.s., {\tt1e-03} & 3   & 3 (74), 3  & {\tt2.05e-05} & 83\\
ROS2, $\tau = T/320$              & 320 &  320 (640), 640 & {\tt2.11e-05}  & 2013 \\
ROS2, $\tau = T/640$              & 640 & 640 (1280), 1280 & {\tt1.66e-05}  & 4368
\\\hline
\multicolumn{5}{c}{\added{grid $60\times 60\times 60$, $T=1\cdot 10^{-3}$, $T\|A(y(T))\|_1\approx\texttt{1.50e+05}$}}\\
\added{nonlin.EBK($m=5$), 10 t.s., {\tt1e-04}}
                             & 30 & 30 (671), 30 & {\tt1.18e-05} & 800
\\
\added{ROS2, $\tau=T/3200$}
                           & 3200 & 3200 (6400), 3200 & {\tt2.18e-06} & 20017
\\\hline
\end{tabular}
\end{table}
%
%
%

Results presented in Table~\ref{t:3D_Bratu} for the nonlinear EBK solver 
are obtained with the block size $m=4$ and $m=5$.  The value $m=4$
is marginally acceptable for the accuracy point of view but does
yield an acceptable error values.  Therefore, the plots in 
Figure~\ref{f:ebk_Bratu}, made to get an insight 
in the error convergence, are produced with a larger block
size \replaced{$m=32$ and $n_s=300$ samples}{$m=12$}.
\added{A convergence stagnation is clearly visible in both plots.
Although we do not have a precise explanation for the stagnation,
we note that it disappears if we switch off the strong 
anisotropy in~\eqref{anis}, taking $A$ to be a Laplacian.
The accuracy level at which stagnation occurs, appears to be related
to the condition number of the Jacobian of the ODE system being solved.
The nonlinear EBK solver relies on repeatedly made low-rank
approximations which are instrumental for the inner block
Krylov subspace solver.  Another source of possible inaccuracies 
is linear system solution within the block Krylov subspace solver
where linear systems with matrices $I+\gamma A_k$ are solved for 
$\gamma>0$ orders of magnitude larger than typical time 
step sizes in implicit time integrators.
Nevertheless, the convergence stagnation is observed at an accuracy
level $\approx 10^{-8}$ which seems to be by far acceptable for PDE solvers.}

\begin{figure}
\centerline{\includegraphics[width=0.48\linewidth]{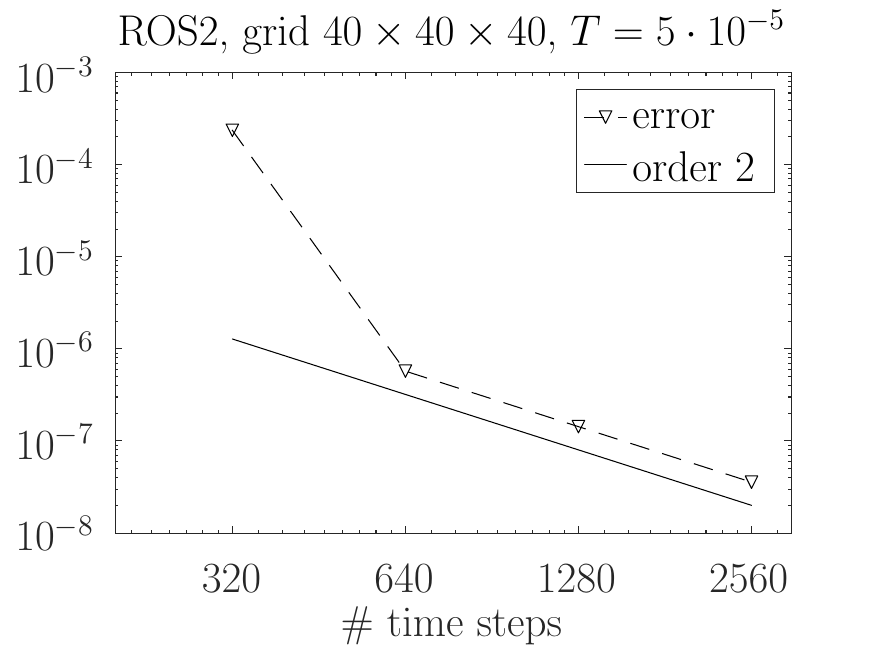}}
\caption{The Bratu 3D test. 
Relative error of the ROS2 method \added{at $t=T$}, computed
as given in~\eqref{err_reached}, versus the time step \replaced{number}{size}.}
\label{f:ros2conv}  
\end{figure}

Finally, we briefly comment on convergence of the ROS2 method.
Since its convergence order can not be seen from error values
shown in Table~\ref{t:3D_Bratu}, we give an error convergence plot 
in Figure~\ref{f:ros2conv}, where relative error values, computed
as given in~\eqref{err_reached}, are plotted versus the 
\replaced{total number of time steps}{corresponding time step sizes}. 
The second order convergence is clearly
seen.

\begin{table}
\caption{\added{Number of nonlinear EBK iterations needed to converge
to relative tolerance $\texttt{tol}=10^{-3}$,
cf.~\eqref{rel_stop}, for different grid sizes and $T$}}
\label{t:Bratu_conv}
\begin{tabular}{cccc}
\hline\hline
      & $20\times 20\times 20$ & $40\times 40\times 40$ & $60\times 60\times 60$ 
\\\hline
$T=2.5\cdot 10^{-5}$ & 2 & 2 & 2 \\
$T=5\cdot 10^{-5}$   & 2 & 2 & 3 \\
$T=1\cdot 10^{-4}$   & 3 & 3 & 3 \\\hline
\end{tabular}
\end{table}

\subsection{3D nonlinear heat conduction}
\label{s:num_heat}
\added{In the previous tests the ODE right hand side function
has a linear part (namely, $-\Asymm y(t)$ in the first test and
$-A y(t)$ in the second one, see~\eqref{1D_discr} 
and~\eqref{3D_discr}, respectively).
To examine whether our nonlinear 
EBK solver is able to handle nonlinear problems with no 
linear part, we now consider the following nonlinear heat
conduction problem, described in~\cite{Zhukov_ea2023}:
find unknown $u=u(x,y,z,t)$ which satisfies
\begin{equation}
\label{heat}
\begin{gathered}
u_t = (k^{(x)}(u)u_x)_x + (k^{(y)}(u)u_y)_y + (k^{(z)}(u)u_z)_z,
\\
(x,y,z)\in[0,1]^3, \quad 0\leqs t\leqs T,
\\
u\bigl|_{x=0} = u\bigl|_{x=1},\quad
u\bigl|_{y=0} = 900,\quad
u\bigl|_{y=1} = 300,\quad
u_z\bigl|_{z=0} = u_z\bigl|_{z=1} = 0,\\
u(x,y,z,0) = 1800\, e^{-60((x-0.5)^2 + (y-0.5)^2 + (z-0.5)^2) },
\\
k^{(x)}(u) = u / 300,
\quad
k^{(y)}(u) = k^{(z)}(u) = k^{(x)}(u)/10.
\end{gathered}
\end{equation}
We discretize this initial boundary value problem on
a regular grid of $n_x\times n_y\times n_z$ nodes
$(x_i,y_j,z_k)$, 
$x_i=i h_x$, $i=1,\dots,n_x$,
$y_j=j h_y$, $j=1,\dots,n_y$,
$z_k=(k-1/2) h_z$, $k=1,\dots,n_z$,
with grid widths $h_{x,y}=1/(n_{x,y}+1)$
and $h_z=1/n_z$.
The shifted grid nodes are taken for the $z$ direction
to accommodate the Neumann boundary conditions.
The standard second order finite difference approximation
is employed, which reads, for the $x$ derivative,
$$
(k^{(x)}(u)u_x)_x \bigl|_{i,j,k} \approx 
\frac{
k^{(x)}(u_{i+1/2,j,k})(u_{i+1,j,k}-u_{i,j,k}) - 
k^{(x)}(u_{i-1/2,j,k})(u_{i,j,k}-u_{i-1,j,k})
}{h_x^2},
$$ 
where $u_{i\pm 1/2,j,k} = (u_{i\pm 1,j,k} + u_{i,j,k})/2$.
The same finite difference approximations are applied
in $y$ and $z$ directions.
This space discretization of the initial boundary value
problem~\eqref{heat} yields an IVP
\begin{equation}
\label{heat_discr}
y'(t) = - A(y(t)) y(t) + g,\quad y(0)=v,  
\end{equation}
where the entries of the vector function $y(t)$ contain
the numerical solution values
and the constant in time vector $g$ contains
the inhomogeneous boundary condition contributions.
The waveform relaxation iteration is applied with
the splitting~\eqref{split_f} where
we set, for $\bar{y}_k=y_k(T)$, 
$$
A_k = A(\bar{y}_k), \quad
f_k(y) = A(\bar{y}_k) y - A(y) y.
$$ 
In this relation, it is not difficult to recognize 
the same splitting as is used for the
Burgers equation test (with $\Asymm=0$ and $\Askew=A$).}

\added{The EBK inner-outer iterations are applied exactly in the 
same setting as in the previous test, with the
same stopping criteria,
see~\eqref{rel_stop}, \eqref{lin_s2} and the same number
of samples $n_s$.  As in the previous test,
we compare the EBK solver with the ROS2 integrator.
The results of the test runs are
presented in Table~\ref{t:3D_heat}.
As we see, on the grid $40\times 40\times 40$ the two
time integration schemes perform equally well
in terms of the CPU time.  Note that the EBK solver
requires much fewer LU~factorizations.
On the finer $60\times 60\times 60$
the EBK solver clearly outperforms ROS2,
both in terms of the CPU time and required LU~factorizations.
When switching to the finer grid, to keep convergence in
the EBK scheme, we have to decrease the time interval $T$
from $T=0.01$ to $T=0.005$.
Hence, we observe a deterioration of the EBK convergence
in this test.  There is also a moderate increase in the 
required number of nonlinear iterations: with
$T=0.005$ EBK needs in total~73 and 94~iterations on the two 
grids, respectively.}

\begin{table}
\caption{\added{The 3D nonlinear heat conduction test.  Attained error and 
computational work for our nonlinear EBK method and 
the ROS2 solver. 
For the EBK method the work is reported as the nonlinear iteration number,
number of LU~factorizations, their applications and matvecs.
For the ROS2 solver the work is shown as 
as the number of time steps, LU~factorizations, their applications and fevals.}}
\label{t:3D_heat}  

\renewcommand{\arraystretch}{1.1} 
\begin{tabular}{ccccc}
\hline\hline
 method, time steps (t.s.) 
                     & iter./ & LUs (LUs applic.), & relative & CPU \\ 
 tolerance / $\tau$  & steps  & matvecs/fevals     & error    & time, s\\
\hline
\multicolumn{5}{c}{\added{grid $40\times 40\times 40$, $T=0.1$, $T\|A(y(T))\|_1\approx\texttt{2.19e+03}$}}\\
nonlin.EBK($m=6$), 10~t.s., {\tt1e-02} 
                     &   52  &  52 (1223), 52 & {\tt1.30e-04}  & 387 \\
nonlin.EBK($m=10$), 10~t.s., {\tt1e-02} 
                     &   52  &  52 (1936), 52 & {\tt5.91e-05}  & 448 \\
ROS2, $\tau=T/250$   &  250  & 250 (500)      & {\tt7.26e-05}  & 347 \\
ROS2, $\tau=T/500$   &  500  & 500 (1000)     & {\tt2.02e-05}  & 821 \\
\hline
\multicolumn{5}{c}{\added{grid $60\times 60\times 60$, $T=0.1$, $T\|A(y(T))\|_1\approx\texttt{4.96e+04}$}}\\
nonlin.EBK($m=6$), 20~t.s., {\tt1e-02} 
                     &   94  &  94 (2068), 94 & {\tt4.41e-05}  & 8229 \\                
ROS2, $\tau=T/500$   &  500   & 500 (1000)    & {\tt3.99e-05}  & 19578 \\\hline
\end{tabular}
\end{table}

\section{Conclusions}
\label{s:concl}

In this paper nonlinear waveform relaxation iterations and their block Krylov 
subspace implementation (the nonlinear EBK method)
are examined theoretically and practically.
Theoretically we have shown that convergence of the considered 
nonlinear waveform relaxation iteration
is determined by the Lipschitz constant $L$ of the nonlinear term $f_k$
in~\eqref{split_f}, the logarithmic norm bound $\omega$ of the linear term 
matrix $A_k$ and the time interval length $T$, see relation~\eqref{conv}.
We have also established a superlinear convergence of our nonlinear iteration, 
which is known to hold in the linear case, see Proposition~\ref{prop2}.
Furthermore, we have linked the residual convergence in the nonlinear
iteration to the error convergence, see relation~\eqref{r2e:eq}.
Finally, it is shown that inexact solution of 
the linear initial-value problem arising at each waveform relaxation
iteration leads to a process whose convergence is equivalent to 
convergence of the exact waveform relaxation iteration with an increased
Lipschitz constant $L+l$, with $l$ being the inexactness tolerance.   

Practical performance of nonlinear waveform relaxation heavily depends on
the ability to solve the arising linear initial-value problem at each 
relaxation iteration quickly and to handle the iterative across-time solutions
efficiently.  As we show here, both these goals can be successfully
reached by employing the exponential block Krylov subspace method
EBK~\cite{Botchev2013}.
In the presented experiments the EBK method is used in the shift-and-invert
(SAI) mode, i.e., it has been assumed that the SAI linear systems (similar to
those arising in implicit time stepping) can be solved efficiently.

Comparisons with implicit time integration methods demonstrated 
superiority and potential of our EBK-SAI based nonlinear waveform
relaxation iterations.  
The presented nonlinear residual stopping criterion has been proven to
be a reliable way to control convergence and to estimate
accuracy.  
Moreover, as \added{the first two} numerical tests demonstrate, 
the grid independent 
convergence of linear SAI Krylov subspace exponential time integrators
(see \cite{MoretNovati04,EshofHochbruck06,Gri12})
\replaced{can be}{is} inherited by the nonlinear waveform relaxation iterations.

\added{As numerical tests and our experience suggest,
the presented method has its limitations.
First, it is not intended for strongly nonlinear problems with
rapidly changing solutions.  The method clearly profits if the right hand 
side function of the ODE system contains a linear part.
Second, the linear systems arising in the presented nonlinear EBK method 
have a structure similar to linear systems occuring in implicit
time integrators with significantly increased time step sizes.
Therefore, the overall efficiency of the method for a particular
application heavily depends on availability of efficient linear solvers.
Third, the nonlinear EBK method is aimed at solving initial-value
problems with a moderate accuracy (of order $10^{-4}$ or $10^{-5}$), 
which is usually enough for solving PDE problems. 
The accuracy of the presented method can be limited by the built-in
low rank approximations and solving large ill-conditioned linear systems.
Finally,} 
some tuning (such as a proper choice of the block size 
$m$ and the time interval length\deleted{ $T$}) is required to apply the proposed approach
successfully. \added{Some of} these issues can be addressed in more details
in a future. 

\added{Nevertheless, as results in this paper show, the presented framework
seems to be promising, has a potential and is of interest for the 
numerical community.}
Recently, an alternative approach based on exponential Krylov subspace
integrators has been presented in~\cite{BotchevZhukov2023,BotchevZhukov2024} 
for nonlinear heat equations.  In this problem a rapid change in the solution
seems to exclude using waveform relaxation with block Krylov subspaces
as presented here.  
A further research should be done to get an insight on when to switch between
the block Krylov subspace approach discussed here and the approach of~\cite{BotchevZhukov2023}.
Finally, it seems to be useful to gain more insight how these two approaches
compare to other exponential time integrators for nonlinear 
initial-value 
problems~\cite{HochbruckOstermann2010,HippHochbruckOstermann14,HansenOstermann2016}. 
 
\subsubsection*{Acknowledgments}
The author thanks the anonymous reviewer for careful reading and suggesting
several valuable improvements.
We acknowledge the use of the hybrid supercomputer K-100 facilities
installed in the Supercomputer center of collective usage 
at Keldysh Institute of Applied Mathematics of Russian Academy of Sciences.

\bibliography{matfun,my_bib}

\end{document}